\documentclass[12pt, a4paper]{article}
\usepackage[latin2]{inputenc}
\usepackage{amssymb}
\usepackage{paralist}
\usepackage{amsmath, calc}
\usepackage{amsfonts}
\usepackage{amsthm}
\usepackage{fullpage}
\usepackage{enumerate}
\usepackage{paralist}
\usepackage{comment}

\newtheorem{theorem}{Theorem}[section] %
\newtheorem{lemma}[theorem]{Lemma} %
\newtheorem{corollary}[theorem]{Corollary} %
\newtheorem{proposition}[theorem]{Proposition} %
\newtheorem{problem}[theorem]{Problem} %
\newtheorem{claim}{Claim} %

\usepackage{cite}

\title{Identical representation functions of linear forms}

\begin{document}

\author{
S\'andor Z. Kiss \thanks{Department of Algebra and Geometry, Institute of Mathematics, Budapest University of Technology and Economics, M\H{u}egyetem rkp. 3., H-1111, Budapest, Hungary. \newline
HUN-REN-BME Stochastics Research Group, M\H{u}egyetem rkp. 3., H-1111 Budapest, Hungary.
Email: kiss.sandor@ttk.bme.hu.
This author was supported by the NKFIH Grants No. K129335, K146387, KKP 144059.}, Csaba
S\'andor \thanks{Department of Stochastics, Institute of Mathematics, Budapest University of
Technology and Economics, M\H{u}egyetem rkp. 3., H-1111, Budapest, Hungary. \newline HUN-REN Alfr\'ed R\'enyi Institute of Mathematics, Re\'altanoda utca 13--15., H-1053 Budapest, Hungary.
MTA--HUN-REN RI Lend\"ulet ``Momentum'' Arithmetic Combinatorics Research Group, Re\'altanoda utca 13--15., H-1053 Budapest, Hungary. Email: sandor.csaba@ttk.bme.hu.
This author was supported by the NKFIH Grants No. K129335, K146387, KKP 144059.}
}

\date{}

\maketitle

\begin{abstract}
\noindent 
For a set of natural numbers $A$, let $R_{A}(n)$ be the number of representations of a natural number $n$ as the sum of two terms from $A$. Many years ago, Nathanson studied the conditions for the set $A$ and $B$ of natural numbers that are needed to guarantee that $R_{A}(n) = R_{B}(n)$ for every positive integer $n$. In the last decades, similar questions have been studied by many authors. In this paper, we extend Nathanson's result to representation functions associated to linear forms and we study related problems.


\end{abstract}

\section{Introduction}

Let $\mathbb{N}$ be the set of nonnegative integers, and let $\mathbb{Z}^{+}$ be the set of positive integers. Let $A\subseteq \mathbb {N}$, $A = \{a_{1}, a_{2}, \dots{}\}$ $(a_{1} < a_{2} < \dots{})$ and $\underline{\alpha} = (\alpha_{1}, \alpha_{2}, \dots{}, \alpha_{s})$ be a vector, where $1 \le \alpha_{1} \le \alpha_{2} \le \dots{} \le \alpha_{s}$ are positive integers. 
For every natural number $n$, the representation function associated to the vector $\underline{\alpha}$ is defined as
\[
R_{A,\underline{\alpha}}(n) = |\{(j_{1}, j_{2}, \dots{} ,j_{s}): j_{i}\in\mathbb{Z}^{+}, \sum_{i=1}^{s}\alpha_{i}a_{j_{i}} = n\}|.
\]
In this paper we study the following problem: 
\begin{problem}\label{prob_1}
    What conditions on $A,B\in \mathbb{N}$, $|A|, |B| \ge 2$ and $\underline{\alpha} = (\alpha_{1}, \alpha_{2}, \dots{}, \alpha_{s})$, $1 \le \alpha_{1} \le \alpha_{2} \le \dots{} \le \alpha_{s}$,  $\underline{\smash{\beta}} = (\beta_{1}, \beta_{2}, \dots{}, \beta_{t})$, $1 \le \beta_{1} \le \beta_{2} \le \dots{} \le \beta_{t}$ are needed to ensure that $R_{A,\underline{\alpha}}(n) = R_{B,\underline{\smash{\beta}}}(n)$ for every nonnegative integer $n$?
\end{problem}

In the last few decades, similar questions have been investigated by many authors [1-6], [8-14]. The first result in this direction was proved by Nathanson in 1978 [9].
If $\underline{\alpha} = (\underbrace{1,1,\dots{} ,1}_{k})$, then for a nonnegative integer $n$, we write $R_{A,\underline{\alpha}}(n) = R_{A,k}(n)$. 

\begin{theorem}[Nathanson, 1978]
For any $A,B\in \mathbb{N}$, one has $R_{A,k}(n) = R_{B,k}(n)$ for every nonnegative integer $n$ if and only if $A = B$.
\end{theorem}

For an integer $m$, a positive integer $k$ and $A,B\subseteq \mathbb{N}$, we define the sets
\[
A + m := \{a_{i} + m: a_{i}\in A\},
\]
\[
kA := \{ka_{i}: a_{i}\in A\},
\]
and the multiset
\[
A + B := \{a_{i} + b_{j}: a_{i}\in A, b_{j}\in B\}. 
\]
Moreover, define the multiset of nonnegative integers $\underline{\alpha}A$ by $\underline{\alpha}A = \alpha_{1}A + \alpha_{2}A + \dots{} + \alpha_{s}A$.
Let
\[
\mathbb{Z}^{+}_{m} = \{\underline{\alpha}: \underline{\alpha} = (\alpha_{1}, \alpha_{2}, \dots{}, \alpha_{s}), 1 \le \alpha_{1} \le \alpha_{2} \le \dots{} \le \alpha_{s}, \alpha_{i}\in \mathbb{Z}^{+}\},
\]
where $s$ is not fixed.
We will show that if for $A,B\subseteq \mathbb{N}$, $A = \{a_{1}, a_{2}, \dots{}\}$ $(a_{1} < a_{2} < \dots{})$,  $B = \{b_{1}, b_{2}, \dots{}\}$ $(b_{1} < b_{2} < \dots{})$,  $\underline{\alpha}, \underline{\smash{\beta}}\in \mathbb{Z}^{+}_{m}$ one has $R_{A,\underline{\alpha}}(n) = R_{B,\underline{\smash{\beta}}}(n)$ for every nonnegative integer $n$, then $R_{A-a_{1},\underline{\alpha}}(n) = R_{B-b_{1},\underline{\smash{\beta}}}(n)$ holds for every nonnegative integer $n$ as well. So we can assume that $0 \in A\cap B$. If $A,B\subseteq \mathbb{N}$, $0 \in A\cap B$ and $\underline{\alpha}\in \mathbb{Z}^{+}_{m}$, then the multiplicity of $0$ is $1$ in the multiset $\underline{\alpha}A$. 

For a finite or infinite multiset $M = \{m_{1}, m_{2}, \dots{}\}$ $(m_{1} \le m_{2} \le \dots{})$ of nonnegative integers, the characteristic function of the multiset $M$ is
\[
\chi_{M}(n) = |\{i: m_{i} = n\}|
\]
for every $n\in \mathbb{N}$. So for every $A\subseteq \mathbb{N}$, $0\in A$ and $\underline{\alpha}\in \mathbb{Z}^{+}_{m}$, we have $\chi_{\underline{\alpha}A}(0) = 1$ and $\chi_{\underline{\alpha}A}(n) < \infty$ for every positive integer $n$. 

We define the set of multisets by
\[
\mathbb{N}_{m} = \{M: \chi_{M}(0) = 1, \chi_{M}(n) < \infty \textnormal{ for every } n\in \mathbb{Z}^{+}\}.
\]
In this paper we focus on the multisets of $\mathbb{N}_{m}$. 
Let $\underline{\alpha} = (\alpha_{1}, \alpha_{2}, \dots{} ,\alpha_{s})$ be a finite vector, where $1 \le \alpha_{1} \le \alpha_{2} \le \dots{} \le \alpha_{s}$ are positive integers.  The characteristic function of $\underline{\alpha}$ is defined by
\[
\chi_{\underline{\alpha}}(n) = |\{i: \alpha_{i} = n\}|.
\]
%Let
%\[
%\mathbb{Z}^{+}_{m} = \{\underline{\alpha}: \chi_{\underline{\alpha}}(n) < \infty \textnormal{ for %every } n\in \mathbb{Z}^{+}\}.
%\]
%For a positive integer $n$ and $A,B\in \mathbb{N}_{m}$, we define the multisets
%\[
%A + n := \{a_{i} + n: a_{i}\in A\},
%\]
%\[
%nA := \{na_{i}: a_{i}\in A\},
%\]
%\[
%A + B := \{a_{i} + b_{j}: a_{i}\in A, b_{j}\in B\}. 
%\]
For multisets of nonnegative integers $A, B$, one can define the multisets $A + B$, $A + m$ and $kA$ in the same way as for the sets of nonnegative integers. 
If $\underline{\alpha} = (\alpha_{1}, \alpha_{2}, \dots{} ,\alpha_{s})$ and $A\in \mathbb{N}_{m}$, then for $\underline{\alpha}A = \alpha_{1}A + \alpha_{2}A + \dots{} 
 + \alpha_{s}A$, we have $\underline{\alpha}A\in \mathbb{N}_{m}$. If $\underline{\alpha}A$ is a set of nonnegative integers, then we write $\underline{\alpha}A \subseteq \mathbb{N}$.
For a positive integer $n$ and a multiset $A \in \mathbb{N}_{m}$, let us define the multiset
\[
A \cap [n]_{m} = \{a_{i}\in A: a_{i}\le n\}.
\]
Similarly, for any vector $\underline{\alpha} = (\alpha_{1}, \alpha_{2}, \dots{} ,\alpha_{s})$, $1 \le \alpha_{1} \le \alpha_{2} \le \dots{} \le \alpha_{s}$, let
\[
\underline{\alpha}\cap [n]_{m} = (\alpha_{1}, \alpha_{2}, \dots{}, \alpha_{t}),
\]
where $\alpha_{t} \le n$ and $\alpha_{t+1} > n$ or $\alpha_{t+1}$ does not exist.
We define the limit of the multisets in the following way. For a sequence of multisets $M^{(1)}, M^{(2)}, \dots{}$, we write
\[
\lim_{k\rightarrow \infty}M^{(k)} = M
\]
if for every positive integer $n$, there exists a positive integer $N$ such that for every $k \ge N$, then $M^{(k)} \cap [n]_{m} = M\cap [n]_{m}$.
Moreover, for a sequence of vectors $\underline{v}_{1}, \underline{v}_{2}, \dots{},$ we write
\[
\lim_{k\rightarrow \infty}\underline{v}_{k} = \underline{v}
\]
if for every positive integer $n$, there exists a positive integer $N$ such that for every $k \ge N$, then $\underline{v}_{k} \cap [n]_{m} = \underline{v}\cap [n]_{m}$. 
%If $A\in \mathbb{N}_{m}$, $\underline{\alpha}\in \mathbb{Z}^{+}_{m}$ is infinite, then let
%\[
%\underline{\alpha}A = \lim_{n\rightarrow \infty}(\underline{\alpha}\cap [n]_{m})A.
%\]
%Obviously, $\underline{\alpha}A \in \mathbb{N}_{m}$.
If $A$ is a multiset of nonnegative integers with $\chi_{A}(n) < \infty$ for every nonnegative integer $n$, $\underline{\alpha}\in \mathbb{Z}^{+}_{m}$, $(\alpha_{1}, \alpha_{2}, \dots{}, \alpha_{s})$ then, the representation function is defined as
\[
R_{A,\underline{\alpha}}(n) = |\{(j_{1}, j_{2}, \dots{} ,j_{s}): j_{i}\in\mathbb{Z}^{+}, \sum_{i=1}^{s}\alpha_{i}a_{j_{i}} = n\}|.
\]
Evidently, for every natural number $n$, one has
\begin{equation}\label{eq:1}
R_{A,\underline{\alpha}}(n) = \chi_{\underline{\alpha}A}(n) = R_{\underline{\alpha}A, (1)}(n). 
\end{equation}
If $A\in \mathbb{N}_{m}$, $\underline{\alpha}\in\mathbb{Z}_{m}^{+}$, then $R_{A,\underline{\alpha}}(0) = 1$ and $R_{A,\underline{\alpha}}(n) < \infty$ 
for every positive integer $n$, that is $\underline{\alpha}A\in \mathbb{N}_{m}$.
If $A,B\in \mathbb{N}_{m}$, $\underline{\alpha},\underline{\smash{\beta}}\in \mathbb{Z}^{+}_{m}$ and $\underline{\alpha}A = \underline{\smash{\beta}}B$, that is
$R_{A,\underline{\alpha}}(n) = R_{B,\underline{\smash{\beta}}}(n)$ for every $n\in \mathbb{N}$, then we simply write
\begin{equation}\label{eq:2}
R_{A,\underline{\alpha}} = R_{B,\underline{\smash{\beta}}}. 
\end{equation}
In this paper we study the following general problem: 
\begin{problem}\label{prob_1}
    What conditions on $A,B\in \mathbb{N}_{m}$, $|A|, |B| \ge 2$ and $\underline{\alpha},\underline{\smash{\beta}}\in \mathbb{Z}^{+}_{m}$ are needed to ensure that $\underline{\alpha}A = \underline{\smash{\beta}}B$, that is $R_{A,\underline{\alpha}} = R_{B,\underline{\smash{\beta}}}$?
\end{problem}

%It is clear that if $\underline{\alpha}$ or $ \underline{\smash{\beta}}$ are infinite in %\eqref{eq:2} and $0\notin A$, then $0\notin B$ and $R_{A,\underline{\alpha}}(n) = %R_{B,\underline{\smash{\beta}}}(n) = 0$ for every $n \in \mathbb{N}$, thus we need $0 \in A\cap B$. 

%If $\underline{\alpha}$ and 
%$\underline{\smash{\beta}}$ are finite, then 

It may happen that $0\notin A\cap B$.
Let us suppose that $R_{A,\underline{\alpha}} = R_{B,\underline{\smash{\beta}}}$.
The next theorem describes all the pairs $(u,v)$ with $R_{A+u,\underline{\alpha}} = R_{B+v,\underline{\beta}}$.

\begin{theorem}\label{prop_1}
    Let $A = \{a_{1}, a_{2}, \dots{} \}$ $(0 \le a_{1} \le a_{2} \le \dots{})$, 
    $B = \{b_{1}, b_{2}, \dots{}\}$ $(0\le b_{1}\le b_{2}  \le \dots{})$ and
    let $\underline{\alpha} = (\alpha_{1}, \alpha_{2}, \dots{} ,\alpha_{s})$,  $\underline{\smash{\beta}} = (\beta_{1}, \beta_{2}, \dots{} ,\beta_{t})$ and $A,B$ be multisets of nonnegative integers such that $R_{A,\underline{\alpha}}(n) = R_{B,\underline{\smash{\beta}}}(n)$ holds for every nonnegative integer $n$. 
    \begin{itemize}
    \item[(1)] If $u \ge -a_{1}$, $v \ge -b_{1}$, $u,v\in \mathbb{Z}$, then $R_{A+u,\underline{\alpha}}(n) = R_{B+v,\underline{\beta}}(n)$ holds for every nonnegative integer $n$ if and only if there exists a $w\in \mathbb{Z}$ with
\[
u = w\frac{\beta_{1} + \beta_{2} + \dots{} + \beta_{t}}{\textnormal{gcd}(\alpha_{1} + \alpha_{2} + \dots{} + \alpha_{s}, \beta_{1} + \beta_{2} + \dots{} + \beta_{t})},
\]
\[
v = w\frac{\alpha_{1} + \alpha_{2} + \dots{} + \alpha_{s}}{\textnormal{gcd}(\alpha_{1} + \alpha_{2} + \dots{} + \alpha_{s}, \beta_{1} + \beta_{2} + \dots{} + \beta_{t})}.
\]
\item[(2)] The equation $R_{A-a_{1},\underline{\alpha}}(n) = R_{B-b_{1},\underline{\smash{\beta}}}(n)$ holds for every nonnegative integer $n$. 
\end{itemize}
\end{theorem}
As a corollary, we get that
\begin{corollary}\label{cor_1} Let $A,B\in \mathbb{N}_{m}$ and let $\underline{\alpha} = (\alpha_{1}, \alpha_{2}, \dots{} ,\alpha_{s})$, $\underline{\smash{\beta}} = (\beta_{1}, \beta_{2}, \dots{} ,\beta_{t})$, $\underline{\alpha}, \underline{\smash{\beta}}\in \mathbb{Z}^{+}_{m}$ such that \eqref{eq:2} holds. Then 
    $R_{A+u,\underline{\alpha}} = R_{B+v,\underline{\smash{\beta}}}$ if and only if there exists a $w\in \mathbb{N}$ such that 
\[
u = w\frac{\beta_{1} + \beta_{2} + \dots{} + \beta_{t}}{\textnormal{gcd}(\alpha_{1} + \alpha_{2} + \dots{} + \alpha_{s}, \beta_{1} + \beta_{2} + \dots{} + \beta_{t})},
\]
\[
v = w\frac{\alpha_{1} + \alpha_{2} + \dots{} + \alpha_{s}}{\textnormal{gcd}(\alpha_{1} + \alpha_{2} + \dots{} + \alpha_{s}, \beta_{1} + \beta_{2} + \dots{} + \beta_{t})}.
\]
\end{corollary}
We say $R_{A,\underline{\alpha}} = R_{B,\underline{\beta}}$ is reduced, if $0\in A\cap B$. A solution of \eqref{eq:2} is called trivial if $(A,\underline{\alpha}) = (B,\underline{\smash{\beta}})$. Obviously, it is enough to consider the reduced non-trivial solutions of \eqref{eq:2}.

%In the last few decades, similar questions have been investigated by many authors [1-6], [8-14]. %The first result in this direction was proved by Nathanson in 1978 [9].
%If $\underline{\alpha} = (\underbrace{1,1,\dots{} ,1}_{k})$, then we write %$R_{A,\underline{\alpha}} = R_{A,k}$. 

%\begin{theorem}[Nathanson, 1978]
%For any $A,B\in \mathbb{N}$, one has $R_{A,k} = R_{B,k}$ if and only if $A = B$.
%\end{theorem}

In this paper, we extend Nathanson's theorem to linear forms and multisets. 

\begin{theorem}\label{theorem_1}
    Let $A,B\in \mathbb{N}_{m}$ with $|A|, |B|\ge 2$, $\underline{\alpha}, \underline{\smash{\beta}}\in \mathbb{Z}^{+}_{m}$. Then we have
    \begin{itemize}
        \item[(i)] $\underline{\alpha}A = \underline{\alpha}B$ if and only if $A = B$.
         \item[(ii)] $\underline{\alpha}A = \underline{\smash{\beta}}A$
    if and only if $\underline{\alpha} = \underline{\smash{\beta}}$. 
    \end{itemize}
\end{theorem}
We assumed that $|A|, |B|\ge 2$ because if $|A| = |B| = 1$, then $A = B = \{0\}$ and $\underline{\alpha}A = \underline{\smash{\beta}}A$ for every $\underline{\alpha}, \underline{\smash{\beta}}\in \mathbb{Z}^{+}_{m}$.
\begin{corollary}\label{cor_2}
    Let $A,B$ be sets of nonnegative integers with $|A|, |B|\ge 2$ and let $\underline{\alpha}, \underline{\smash{\beta}}\in \mathbb{Z}^{+}_{m}$ be finite vectors, where $\underline{\alpha} = (\alpha_{1}, \dots{} ,\alpha_{s})$ and $\underline{\smash{\beta}} = (\beta_{1}, \dots{} ,\beta_{t})$. Then we have
    \begin{itemize}
        \item[(i)] $\underline{\alpha}A = \underline{\alpha}B$ if and only if $A = B$.
         \item[(ii)] $\underline{\alpha}A = \underline{\smash{\beta}}A$
    if and only if $\underline{\alpha} = \underline{\smash{\beta}}$. 
    \end{itemize}
\end{corollary}

As a corollary, we immediately get that

\begin{corollary}\label{cor_3}
    If \eqref{eq:2} has a nontrivial solution, $A, B\in \mathbb{N}_{m}$, $\underline{\alpha}, \underline{\smash{\beta}}\in \mathbb{Z}_{m}^{+}$, where $|A|, |B|\ge 2$, then 
$A \neq B$ and $\underline{\alpha} \neq \underline{\smash{\beta}}$.
\end{corollary}

For the sake of simplicity, we denote the ordered Kronecker product of two vectors $\underline{\alpha} = (\alpha_{1}, \alpha_{2}, \dots{} ,\alpha_{s})$ and 
$\underline{\smash{\beta}} = (\beta_{1}, \beta_{2}, \dots{} ,\beta_{t})$ by 
$\underline{\alpha}\,\underline{\smash{\beta}} = (\alpha_{1}\beta_{1}, \dots{}  ,\alpha_{i}\beta_{j},
\dots{} ,\alpha_{s}\beta_{t})$. 

If $\underline{\alpha} = (\alpha_{1}, \alpha_{2}, \dots{} ,\alpha_{s})$, $\underline{\smash{\beta}} = (\beta_{1}, \beta_{2}, \dots{} ,\beta_{t})$ and $A\in \mathbb{N}_{m}$ then
\[
(\underline{\alpha}\,\underline{\smash{\beta}})A = \sum_{(i,j)}(\alpha_{i}\beta_{j})A = \sum_{i=1}^{s}\sum_{j=1}^{t}(\alpha_{i}\beta_{j})A = \sum_{i=1}^{s}\alpha_{i}\left (\sum_{j=1}^{t}\beta_{j}A \right) = \sum_{i=1}^{s}\alpha_{i}(\underline{\smash{\beta}} A) = \underline{\alpha}(\underline{\smash{\beta}}A).
\]
%If $\underline{\alpha}$ or $\underline{\smash{\beta}}$ is infinite, then for every positive %integer $n$, 
%\[
%(\underline{\alpha}\,\underline{\smash{\beta}})A \cap [n]_{m} = ((\underline{\alpha}\cap [n]_{m})
%(\underline{\smash{\beta}}\cap [n]_{m})(A\cap [n]_{m}))\cap [n]_{m}
%\]
%\[
%= ((\underline{\alpha}\cap [n]_{m})(\underline{\smash{\beta}}A\cap [n]_{m}))\cap [n]_{m} = %\underline{\alpha}(\underline{\smash{\beta}}A)\cap [n]_{m}.
%\]
%Since the above equation holds for every natural number $n$, we have %$(\underline{\alpha}\,\underline{\smash{\beta}})A = \underline{\alpha}(\underline{\smash{\beta}}A)$. 
By the definition, for every $\underline{\alpha}, \underline{\smash{\beta}}, \underline{\smash{\gamma}}\in \mathbb{Z}_{m}^{+}$, we have $\underline{\alpha} = \underline{\alpha}(1)$,
$\underline{\alpha}\,\underline{\smash{\beta}} = \underline{\smash{\beta}}\underline{\alpha}$, 
$(\underline{\alpha}\,\underline{\smash{\beta}})\underline{\smash{\gamma}} = \underline{\alpha}(\underline{\smash{\beta}}\,\underline{\smash{\gamma}})$. Now we get a nontrivial solution of \eqref{eq:2} in the following way. Let $\underline{\alpha}_{1}, \underline{\alpha}_{2}, \underline{\smash{\beta}}_{1}, \underline{\smash{\beta}}_{2}, \in \mathbb{Z}^{+}$ with $\underline{\alpha}_{1}\underline{\alpha}_{2} = \underline{\smash{\beta}}_{1}\underline{\smash{\beta}}_{2}$.
Then for $C\in \mathbb{N}_{m}$, we have $(\underline{\alpha}_{1}\underline{\alpha}_{2})C = \underline{\alpha}_{2}(\underline{\alpha}_{1}C)$, that is  
\[
R_{(\underline{\alpha}_{1}\underline{\alpha}_{2})C,(1)} =  R_{\underline{\alpha}_{2}(\underline{\alpha}_{1}C),(1)} = R_{\underline{\alpha}_{1}C, \underline{\alpha}_{2}}
\]
and $(\underline{\smash{\beta}}_{1}\underline{\smash{\beta}}_{2})C = \underline{\smash{\beta}}_{2}(\underline{\smash{\beta}}_{1}C)$, that is
\[
R_{(\underline{\smash{\beta}}_{1}\underline{\smash{\beta}}_{2})C,(1)} =  R_{\underline{\smash{\beta}}_{2}(\underline{\smash{\beta}}_{1}C),(1)} = R_{\underline{\smash{\beta}}_{1}C, \underline{\smash{\beta}}_{2}}.
\]
This allows us to define a nontrivial solution of \eqref{eq:2}. 

\textbf{First construction.} Let us suppose that for
$\underline{\alpha}_{1}, \underline{\alpha}_{2}, \underline{\smash{\beta}}_{1}, \underline{\smash{\beta}}_{2}\in \mathbb{Z}_{m}^{+}$ and $C\in \mathbb{N}_{m}$ we have
$\underline{\alpha}_{1}\underline{\alpha}_{2} = \underline{\smash{\beta}}_{1}\underline{\smash{\beta}}_{2}$, $(\underline{\alpha}_{1}, \underline{\alpha}_{2}) \neq 
(\underline{\smash{\beta}}_{1}, \underline{\smash{\beta}}_{2})$. Then 
$\underline{\alpha}_{2}(\underline{\alpha}_{1}C) = \underline{\smash{\beta}}_{2}(\underline{\smash{\beta}}_{1}C)$, that is
\begin{equation}\label{eq:3}
    R_{\underline{\alpha}_{1}C, \underline{\alpha}_{2}} = R_{\underline{\smash{\beta}}_{1}C, \underline{\smash{\beta}}_{2}}.
\end{equation}
Obviously, as a special case we get that for any $\underline{\alpha}_{1}, \underline{\alpha}_{2}\in \mathbb{Z}^{+}_{m}$, $\underline{\alpha}_{1}\underline{\alpha}_{2} \neq (1)$ and $C \in \mathbb{N}_{m}$ we have $R_{\underline{\alpha}_{1}C, \underline{\alpha}_{2}} = R_{\underline{\alpha}_{2}C, \underline{\alpha}_{1}}$. We have a more general way to get a solution of \eqref{eq:2}. 

%The solutions of the form \eqref{eq:3} are called solutions arising from the first %construction. 
%By the definition of the additive representation function,
%\[
%R_{\underline{\alpha}_{1}A, \underline{\alpha}\,\underline{\alpha}_{2}}(n) = %R_{\underline{\alpha}_{1}\underline{\alpha}_{2}A, \underline{\alpha}}(n) 
%\]
%\[
%R_{\underline{\smash{\beta}}_{1}\underline{\smash{\beta}}_{2}B, %\underline{\smash{\beta}}}(n) =
%R_{\underline{\smash{\beta}}_{1}B, %\underline{\smash{\beta}}\,\underline{\smash{\beta}}_{2}}(n). 
%\]
%It follows that
\begin{proposition}\label{prop_2}
    Let $\underline{\alpha}, \underline{\smash{\beta}}\in \mathbb{Z}^{+}_{m}$, $A,B\in \mathbb{N}_{m}$ such that $\underline{\alpha}A = \underline{\smash{\beta}}B$, that is
    $R_{A,\underline{\alpha}} =  R_{B,\underline{\smash{\beta}}}$. If $\underline{\alpha}_{1}, \underline{\alpha}_{2}, \underline{\smash{\beta}}_{1},  \underline{\smash{\beta}}_{2} \in \mathbb{Z}_{m}^{+}$ with $\underline{\alpha}_{1}\underline{\alpha}_{2} = \underline{\smash{\beta}}_{1}\underline{\smash{\beta}}_{2}$, then $(\underline{\alpha}\,\underline{\alpha}_{2})(\underline{\alpha}_{1}A) = (\underline{\smash{\beta}}\,\underline{\smash{\beta}}_{2})(\underline{\smash{\beta}}_{1}B)$, that is
    $R_{\underline{\alpha}_{1}A, \underline{\alpha}\,\underline{\alpha}_{2}} = R_{\underline{\smash{\beta}}_{1}B, \underline{\smash{\beta}}\,\underline{\smash{\beta}}_{2}}$.
\end{proposition}
Obviously, if $\underline{\alpha} = \underline{\smash{\beta}} = (1)$, $A = B = C$, we get the first construction.
Unfortunately, iterating this proposition to a solution arising from the first construction, we get another solution from the first construction, because if $C\in \mathbb{N}_{m}$, 
$\underline{\alpha}_{1}\underline{\alpha}_{2} = \underline{\smash{\beta}}_{1}\underline{\smash{\beta}}_{2}$, $\underline{\alpha}_{3}\underline{\alpha}_{4} = \underline{\smash{\beta}}_{3}\underline{\smash{\beta}}_{4}$, then $\underline{\alpha}_{2}(\underline{\alpha}_{1}C) = \underline{\smash{\beta}}_{2}(\underline{\smash{\beta}}_{1}C)$
and by Proposition \ref{prop_2}, we have  
$(\underline{\alpha}_{2}\underline{\alpha}_{4})(\underline{\alpha}_{1}\underline{\alpha}_{3}C) = (\underline{\smash{\beta}}_{2}\underline{\smash{\beta}}_{4})(\underline{\smash{\beta}}_{1}\underline{\smash{\beta}}_{3}C)$, where $(\underline{\alpha}_{2}\underline{\alpha}_{4})(\underline{\alpha}_{1}\underline{\alpha}_{3}) = (\underline{\smash{\beta}}_{2}\underline{\smash{\beta}}_{4})(\underline{\smash{\beta}}_{1}\underline{\smash{\beta}}_{3})$.

It is easy to see that the converse of Proposition \ref{prop_2} also holds in the following way. 
If $\underline{\alpha}_{1}\underline{\alpha}_{2} = \underline{\smash{\beta}}_{1}\underline{\smash{\beta}}_{2}$, $A, B \in \mathbb{N}_{m}$, $|A|, |B| \ge 2$ and $(\underline{\alpha}\,\underline{\alpha}_{2})(\underline{\alpha}_{1}A) = (\underline{\smash{\beta}}\,\underline{\smash{\beta}}_{2})(\underline{\smash{\beta}}_{1}B)$, then $(\underline{\alpha}_{1}\underline{\alpha}_{2})(\underline{\alpha} A) = (\underline{\smash{\beta}}_{1}\underline{\smash{\beta}}_{2})(\underline{\smash{\beta}} B)$.
By (i) in Theorem \ref{theorem_1}, we have $\underline{\alpha}A = \underline{\smash{\beta}}B$  and so we get   

\begin{proposition}\label{prop_3}
    If $\underline{\alpha}, \underline{\alpha}_{1}, \underline{\alpha}_{2}, \underline{\smash{\beta}}, \underline{\smash{\beta}}_{1}, \underline{\smash{\beta}}_{2} \in \mathbb{Z}_{m}^{+}$ with $\underline{\alpha}_{1}\underline{\alpha}_{2} = \underline{\smash{\beta}}_{1}\underline{\smash{\beta}}_{2}$, $A, B \in \mathbb{N}_{m}$, $|A|, |B| \ge 2$ and $(\underline{\alpha}\,\underline{\alpha}_{2})(\underline{\alpha}_{1}A) = (\underline{\smash{\beta}}\,\underline{\smash{\beta}}_{2})(\underline{\smash{\beta}}_{1}B)$, then we have $\underline{\alpha}A =  \underline{\smash{\beta}}B$. 
\end{proposition}

It follows from Proposition \ref{prop_3} that if $\underline{\alpha}, \underline{\alpha}_{1}, \underline{\alpha}_{2}, \underline{\smash{\beta}}, \underline{\smash{\beta}}_{1}, \underline{\smash{\beta}}_{2}, \underline{\smash{\gamma}}, \underline{\delta}\in \mathbb{Z}_{m}^{+}$ and 
$A,B,C,D\in \mathbb{N}_{m}$, $|C|, |D|\ge 2$, $\underline{\alpha}A = \underline{\smash{\beta}}B$,
$A = \underline{\alpha}_{1}C$, $\underline{\alpha} = \underline{\smash{\gamma}}\,\underline{\alpha}_{2}$, $B = \underline{\smash{\beta}}_{1}D$, $\underline{\smash{\beta}} =  \underline{\delta}\,\underline{\smash{\beta}}_{2}$ 
and $\underline{\alpha}_{1}\underline{\alpha}_{2} = \underline{\smash{\beta}}_{1}\underline{\smash{\beta}}_{2}$, 
then $(\underline{\smash{\gamma}}\,\underline{\alpha}_{2})(\underline{\alpha}_{1}C) = (\underline{\delta}\,\underline{\smash{\beta}}_{2})(\underline{\smash{\beta}}_{1}D)$. By simplifying with $\underline{\alpha}_{1}, \underline{\alpha}_{2}, \underline{\smash{\beta}}_{1}, \underline{\smash{\beta}}_{2}$, we get $\underline{\smash{\gamma}}C = \underline{\delta}D$. Now we define the primitive solutions of equation \eqref{eq:2}. For $A,B\in \mathbb{N}_{m}$, $\underline{\alpha}, \underline{\smash{\beta}}\in \mathbb{Z}^{+}_{m}$, we say that a solution of \eqref{eq:2} is primitive if the vectors $\underline{\alpha}_{1}, \underline{\alpha}_{2}, \underline{\smash{\beta}}_{1}, \underline{\smash{\beta}}_{2}, \underline{\smash{\gamma}}, \underline{\delta}\in \mathbb{Z}^{+}_{m}$ and multisets $C,D\in \mathbb{N}_{m}$, $|C|, |D|\ge 2$ satisfy $A = \underline{\alpha}_{1}C$, $\underline{\alpha} = \underline{\smash{\gamma}}\,\underline{\alpha}_{2}$, $B = \underline{\smash{\beta}}_{1}D$, $\underline{\smash{\beta}} =  \underline{\delta}\,\underline{\smash{\beta}}_{2}$, and $\underline{\alpha}_{1}\underline{\alpha}_{2} = \underline{\smash{\beta}}_{1}\underline{\smash{\beta}}_{2}$, whenever $\underline{\alpha}_{1} = \underline{\smash{\beta}}_{1} = \underline{\alpha}_{2} = \underline{\smash{\beta}}_{2} = (1)$.  
It is easy to see that 
\begin{claim}\label{claim_1}
Let $A, B \in \mathbb{N}_{m}$, $|A|, |B| \ge 2$. If $\underline{\alpha} A = \underline{\smash{\beta}} B$ is a primitive solution, then this equation cannot arise from the first construction.
\end{claim}
Our next theorem shows that equation \eqref{eq:2} can be simplified so that it gives a primitive solution.

\begin{theorem}\label{theorem_2}
    For any $A,B\in \mathbb{N}_{m}$ with $|A|, |B| \ge 2$, $\underline{\alpha}, \underline{\smash{\beta}}\in \mathbb{Z}^{+}_{m}$ with $\underline{\alpha}A =  \underline{\smash{\beta}}B$, there exist $\underline{\alpha}_{1}, \underline{\alpha}_{2},  \underline{\smash{\beta}}_{1}, \underline{\smash{\beta}}_{2}, \underline{\smash{\gamma}}, \underline{\delta}\in \mathbb{Z}_{m}^{+}$ and 
$C,D\in \mathbb{N}_{m}$, $|C|, |D| \ge 2$ such that $A = \underline{\alpha}_{1}C$, $\underline{\alpha} = \underline{\smash{\gamma}}\,\underline{\alpha}_{2}$, $B = \underline{\smash{\beta}}_{1}D$, $\underline{\smash{\beta}} =  \underline{\delta}\,\underline{\smash{\beta}}_{2}$, $\underline{\alpha}_{1}\underline{\alpha}_{2} = \underline{\smash{\beta}}_{1}\underline{\smash{\beta}}_{2}$ and $\underline{\smash{\gamma}}C = \underline{\delta}D$ is a primitive solution.
\end{theorem}

It follows that it is enough to determine the reduced primitive solutions of $\underline{\alpha} A = \underline{\smash{\beta}} B$ and the vectors  
$\underline{\alpha}_{1}, \underline{\alpha}_{2}, \underline{\smash{\beta}}_{1}, \underline{\smash{\beta}}_{2} \in \mathbb{Z}_{m}^{+}$ with
$\underline{\alpha}_{1}\underline{\alpha}_{2} = \underline{\smash{\beta}}_{1}\underline{\smash{\beta}}_{2}$.

\begin{problem}\label{prob_2}
    Determine all the vectors $\underline{\alpha}_{1}, \underline{\alpha}_{2}, \underline{\smash{\beta}}_{1}, \underline{\smash{\beta}}_{2} \in \mathbb{Z}^{+}_{m}$ with $\underline{\alpha}_{1}\underline{\alpha}_{2} = \underline{\smash{\beta}}_{1}\underline{\smash{\beta}}_{2}$. 
\end{problem}

Since one cannot formulate a unique factorization type result because, for example, $(1,2,4,8,16,32) = (1,2)(1,4,16) = (1,2,4)(1,8)$, which shows that this problem seems to be difficult. 

Now we show that there exist reduced primitive non-trivial solutions of the equation $\underline{\alpha}A = \underline{\smash{\beta}}B$. Let $\lambda \ge 2$ be an integer.
For $N\in \mathbb{Z}^{+}$, let 
\[
S^{(\lambda,N)} = \left \{\sum_{i=0}^{\infty}\delta_{i}\lambda^{iN}: \delta_{i}\in \{0,1,\dots{} ,\lambda-1\}\right \}.
\]
Since there is a unique representation 
\[
n = \sum_{i=0}^{\infty}\alpha_{i}\lambda^{i}, \hspace*{5mm} \alpha_{i}\in \{0,1,\dots{} ,\lambda-1\},
\]
where
\[
n = \sum_{j=0}^{N-1}\lambda^{j}\sum_{i=0}^{\infty}\alpha_{j+iN}\lambda^{iN}
\]
therefore, for every natural number $n$,
\[
R_{S^{(\lambda,N)},(1,\lambda, \lambda^{2}, \dots{},\lambda^{N-1})}(n) = 1.
\]
%If $N$ is infinite, then let $S^{(\lambda,N)} = \{0,1,\dots{} ,\lambda-1\}$ and 
%$(1,\lambda, \lambda^{2}, \dots{} ,\lambda^{N-1}) = (1,\lambda, \lambda^{2}, \dots{})$.
%Then for every natural number $n$,  
%\[
%R_{S^{(\lambda,N)},(1,\lambda, \lambda^{2}, \dots{} ,\lambda^{N-1})}(n) = 1.
%\]
Let $u_{i}, v_{i} \ge 2$ be integers and let $N = u_{1}v_{1}u_{2}v_{2}\cdots{} u_{r}v_{r}$ be a finite product. 
Let
\begin{eqnarray*}
\underline{\lambda}_{\underline{u},\underline{v}}^{(\underline{u})} = 
(1,\lambda, \lambda^{2}, \dots{} ,\lambda^{u_{1}-1})(1,\lambda^{u_{1}v_{1}}, \lambda^{2u_{1}v_{1}}, \dots{} ,\lambda^{(u_{2}-1)u_{1}v_{1}})\\
(1,\lambda^{u_{1}v_{1}u_{2}v_{2}}, \lambda^{2u_{1}v_{1}u_{2}v_{2}}, \dots{} ,\lambda^{(u_{3}-1)u_{1}v_{1}u_{2}v_{2}})\dots{}\\
(1,\lambda^{u_{1}v_{1}u_{2}v_{2}\cdots{} u_{r-1}v_{r-1}}, \lambda^{2u_{1}v_{1}u_{2}v_{2}\cdots{} u_{r-1}v_{r-1}}, \dots{} ,\lambda^{(u_{r}-1)u_{1}v_{1}u_{2}v_{2} \cdots{} u_{r-1}v_{r-1}})
\end{eqnarray*}
\begin{eqnarray*}
\underline{\lambda}_{\underline{u},\underline{v}}^{(\underline{v})} = 
(1,\lambda^{u_{1}}, \lambda^{2u_{1}}, \dots{} ,\lambda^{(v_{1}-1)u_{1}})(1,\lambda^{u_{1}v_{1}u_{2}}, \lambda^{2u_{1}v_{1}u_{2}}, \dots{},
\lambda^{(v_{2}-1)u_{1}v_{1}u_{2}})\\
(1,\lambda^{u_{1}v_{1}u_{2}v_{2}u_{3}}, \lambda^{2u_{1}v_{1}u_{2}v_{2}u_{3}}, \dots{} ,\lambda^{(v_{3}-1)u_{1}v_{1}u_{2}v_{2}u_{3}})\dots{}\\
(1,\lambda^{u_{1}v_{1}u_{2}v_{2}\cdots{}  u_{r-1}v_{r-1}u_{r}}, \lambda^{2u_{1}v_{1}u_{2}v_{2}\cdots{}  u_{r-1}v_{r-1}u_{r}}, \dots{} ,\lambda^{(v_{r}-1)u_{1}v_{1}u_{2}v_{2} \cdots{}  u_{r-1}v_{r-1}u_{r}})
\end{eqnarray*}
Since every $0 \le n \le N - 1$ can be uniquely written in the form
\[
n = \sum_{i=1}^{r}(\delta_{i}u_{1}v_{1}u_{2}v_{2} \cdots{} u_{i-1}v_{i-1} + \delta_{i}^{'}u_{1}v_{1}u_{2}v_{2} \cdots{} u_{i-1}v_{i-1}u_{i}),
\]
where $\delta_{i}\in \{0,1, \dots{} ,u_{i}-1\}$, $\delta_{i}^{'}\in \{0,1, \dots{} ,v_{i}-1\}$ then we have 
\[
\underline{\lambda}_{\underline{u},\underline{v}}^{(\underline{u})}\underline{\lambda}_{\underline{u},\underline{v}}^{(\underline{v})} = (1,\lambda, \lambda^{2}, \dots{} ,\lambda^{N-1}), 
\]
and so for every $n\in \mathbb{N}$ we have
\[
1 = R_{S^{(\lambda,N)},(1,\lambda, \lambda^{2}, \dots{} ,\lambda^{N-1})}(n) = R_{S^{(\lambda,N)},\underline{\lambda}_{\underline{u},\underline{v}}^{(\underline{u})}\underline{\lambda}_{\underline{u},\underline{v}}^{(\underline{v})}}(n) = R_{\underline{\lambda}_{\underline{u},\underline{v}}^{(\underline{u})}S^{(\lambda,N)},\underline{\lambda}_{\underline{u},\underline{v}}^{(\underline{v})}}(n). 
\]
%It follows from \eqref{eq:1} that for every $n\in \mathbb{N}$,
%\[
%R_{\underline{\lambda}_{\underline{u},\underline{v}}^{(\underline{u})}S^{(\lambda,N)},\underline{\lambda}_{\underline{u},\underline{v}}^{(\underline{v})}}(n) = 1,
%\]
%that is
%\[
%\underline{\lambda}^{(\underline{v})}_{\underline{u},\underline{v}}(\underline{\lambda}^{(\underline{u})}_{\underline{u},\underline{v}}S^{(\lambda,N)}) = \mathbb{N}.
%\]
Thus, we get the second type of construction in the following way. 

\textbf{Second construction.} Assume that $\underline{u} = (u_{1},u_{2}, \dots{} ,u_{r})$, $\underline{v} = (v_{1},v_{2}, \dots{} ,v_{r})$, $\underline{u}^{'} = (u_{1}^{'},u_{2}^{'}, \dots{} ,u_{r^{'}}^{'})$, $\underline{v}^{'} = (v_{1}^{'},v_{2}^{'}, \dots{} ,v_{r^{'}}^{'})$, where $u_{i}, u_{i}^{'}, v_{i}, v_{i}^{'} \ge 2$, $\underline{\lambda}_{\underline{u},\underline{v}}^{(\underline{v})} \neq (1)$, $(\lambda^{'})_{\underline{u}^{'},\underline{v}^{'}}^{(\underline{v}^{'})} \neq (1)$. Then for $\lambda, \lambda^{'} \ge 2$ integers, we have
\[
\underline{\lambda}^{(\underline{v})}_{\underline{u},\underline{v}}(\underline{\lambda}^{(\underline{u})}_{\underline{u},\underline{v}}S^{(\lambda,N)}) = 
(\underline{\lambda^{'}})^{(\underline{v}^{'})}_{\underline{u}^{'},\underline{v}^{'}}
((\underline{\lambda^{'}})^{(\underline{u}^{'})}_{\underline{u}^{'},\underline{v}^{'}}S^{(\lambda^{'},N^{'}}) = \mathbb{N},
\]
that is  $R_{\underline{\lambda}_{\underline{u},\underline{v}}^{(\underline{u})}S^{(\lambda,N)},\underline{\lambda}_{\underline{u},\underline{v}}^{(\underline{v})}}(n) = R_{(\underline{\lambda^{'}})^{(\underline{u}^{'})}_{\underline{u}^{'},\underline{v}^{'}}S^{(\lambda^{'},N^{'}}),(\underline{\lambda^{'}})^{(\underline{v}^{'})}_{\underline{u}^{'},\underline{v}^{'}}}(n)$.

\begin{theorem}\label{theorem_3}
    Assume that 
    \begin{equation}\label{eq:4}
R_{\underline{\lambda}^{(\underline{u})}_{\underline{u},\underline{v}}S^{(\lambda,N)},
\underline{\lambda}^{(\underline{v})}_{\underline{u},\underline{v}}} = R_{(\underline{\lambda^{'}})^{(\underline{u}^{'})}_{\underline{u}^{'},\underline{v}^{'}}S^{(\lambda^{'},N^{'})},(\underline{\lambda^{'}})^{(\underline{v}^{'})}_{\underline{u}^{'},\underline{v}^{'}}} = 1
    \end{equation}
a non trivial solution in which $\lambda_{\underline{u},\underline{v}}^{(\underline{v})} \neq (1)$, $\lambda_{\underline{u}^{'},\underline{v}^{'}}^{(\underline{v}^{'})} \neq (1)$. Then the following statements are equivalent:
\begin{itemize}
\item[(i)] $\log \lambda$, $\log \lambda^{'}$ are linearly independent over $\mathbb{Q}$, where $\mathbb{Q}$ denotes the field of rational numbers.
\item[(ii)] \eqref{eq:4} is a primitive solution.
\item[(iii)] \eqref{eq:4} cannot follows from the first type construction.
\end{itemize}
\end{theorem}

\begin{comment}
\begin{theorem}\label{theorem_3}
    Assume that 
    \begin{equation}\label{eq:4}
\underline{\lambda}^{(\underline{u})}_{\underline{u},\underline{v}}(\underline{\lambda}^{(\underline{v})}_{\underline{u},\underline{v}}S^{(\lambda,N)}) = 
\underline{\lambda^{'}}^{(\underline{u}^{'})}_{\underline{u}^{'},\underline{v}^{'}}
(\underline{\lambda^{'}}^{(\underline{v}^{'})}_{\underline{u}^{'},\underline{v}^{'}}S^{(\lambda^{'},N^{'}})) = \mathbb{N}
    \end{equation}
a non trivial solution in which $\lambda_{\underline{u},\underline{v}}^{(\underline{v})} \neq (1)$, $\lambda_{\underline{u}^{'},\underline{v}^{'}}^{(\underline{v}^{'})} \neq (1)$. Then 
\eqref{eq:4} is a primitive solution if and only if
$\log \lambda$, $\log \lambda^{'}$ are linearly independent over $\mathbb{Q}$, where $\mathbb{Q}$ denotes the field of rational numbers.
\end{theorem}
\noindent The following statement follows naturally from Claim 1 and Theorem \ref{theorem_3}.

\begin{corollary}\label{cor_4}
\eqref{eq:4} arise from the first construction if and only if $\log \lambda$, $\log \lambda^{'}$ are linearly dependent over $\mathbb{Q}$.
\end{corollary}
\end{comment}

Obviously,

\begin{corollary}\label{cor_4}
If $\log \lambda$, $\log \lambda^{'}$ are linearly independent over $\mathbb{Q}$, then
\[
R_{\underline{\lambda}^{(\underline{u})}_{\underline{u},\underline{v}}S^{(\lambda,N)},
\underline{\lambda}^{(\underline{v})}_{\underline{u},\underline{v}}} = R_{(\underline{\lambda^{'}})^{(\underline{u}^{'})}_{\underline{u}^{'},\underline{v}^{'}}S^{(\lambda^{'},N^{'})},(\underline{\lambda^{'}})^{(\underline{v}^{'})}_{\underline{u}^{'},\underline{v}^{'}}}
\]    
a non-trivial, reduced and primitive solution of \eqref{eq:2}.
\end{corollary}

\begin{problem}
    Does there exist any other non-trivial reduced and primitive solution of \eqref{eq:2} in $\mathbb{N}_{m}$?
\end{problem}
In this paper we focus on the set $\mathbb{N}_{m}$. One can ask whether $\mathbb{N}_{m}$ is replaced by $\mathbb{N}$. 
By Theorem \ref{prop_1}, we know that for $A,B\subseteq \mathbb{N}$, $A = \{a_{1}, a_{2}, \dots{}\}$ $(a_{1} < a_{2} < \dots{})$ and $B = \{b_{1}, b_{2}, \dots{}\}$ $(b_{1} < b_{2} < \dots{})$,  $\underline{\alpha}, \underline{\smash{\beta}}\in \mathbb{Z}^{+}_{m}$ with $R_{A,\underline{\alpha}} = R_{B,\underline{\smash{\beta}}}$, we have $R_{A-a_{1},\underline{\alpha}} = R_{B-b_{1},\underline{\smash{\beta}}}$. Thus, we may assume that $0 \in A\cap B$. By the first construction we get that if $A,B\subseteq \mathbb{N}$, $0 \in A\cap B$, 
$\underline{\alpha}_{1},\underline{\alpha}_{2}, \underline{\smash{\beta}}_{1}, \underline{\smash{\beta}}_{2}\in \mathbb{Z}^{+}_{m}$ with $\underline{\alpha}_{1}\underline{\alpha}_{2} =\underline{\smash{\beta}}_{1}\underline{\smash{\beta}}_{2}$ and $\underline{\alpha}A =\underline{\smash{\beta}}B$ such that $\underline{\alpha}_{1}A, \underline{\smash{\beta}}_{1}B \subseteq \mathbb{N}$ then $R_{\underline{\alpha}_{1}A,\underline{\alpha}_{2}\underline{\alpha}} = R_{\underline{\smash{\beta}}_{1}B,\underline{\smash{\beta}}_{2}\underline{\smash{\beta}}}$.

\begin{problem}\label{prob_4}
    Does there exist any other non-trivial reduced and primitive solution of \eqref{eq:2} in $\mathbb{N}$ as given in Corollary \ref{cor_4}?
    %Determine all the distinct vectors $\underline{\alpha}, \underline{\smash{\beta}} \in \mathbb{Z}^{+}_{m}$, for which there exist $A, B\subset \mathbb{N}$ with $|A|, |B| \ge 2$ such that $\underline{\alpha}A =\underline{\smash{\beta}}B$. 
\end{problem}

\begin{problem}\label{prob_4}
Does there exist $A, B\subseteq \mathbb{N}$, $0\in A\cap B$, $2 \le |A|, |B| < \infty$ and $\underline{\alpha}, \underline{\smash{\beta}}\in \mathbb{Z}^{+}_{m}$
such that $\underline{\alpha}A = \underline{\smash{\beta}}B$, $(\underline{\alpha}, A)\neq (\underline{\smash{\beta}}, B)$ a primitive solution?
\end{problem}

\begin{comment}
If the subsums of the coordinates of $\underline{\alpha}$ and $\underline{\smash{\beta}}$ are pairwise distinct, then the sets $A = \underline{\smash{\beta}}\{0,1\}$ and $B = \underline{\alpha}\{0,1\}$ are suitable.
Moreover, if $\underline{\alpha} = \underline{\smash{\gamma}}\,\underline{\alpha}^{'}$,
and $\underline{\smash{\beta}} = \underline{\smash{\gamma}}\,\underline{\smash{\beta}}^{'}$ and the subsums of the coordinates of either $\underline{\alpha}^{'}$ or $\underline{\smash{\beta}}^{'}$ are pairwise distinct, then $A = \underline{\smash{\beta}}^{'}\{0,1\}$ and $B = \underline{\alpha}^{'}\{0,1\}$ are suitable.

\begin{problem}\label{prob_5}
Is it true that if $\underline{\alpha}, \underline{\smash{\beta}} \in \mathbb{Z}^{+}_{m}$, $A, B\subset \mathbb{N}$ with $|A|, |B| \ge 2$ such that $\underline{\alpha}A =\underline{\smash{\beta}}B$,  
then there exists a $\underline{\smash{\gamma}}\in \mathbb{Z}_{m}^{+}$ with $\underline{\alpha} = \underline{\smash{\gamma}}\,\underline{\alpha}^{'}$,
$\underline{\smash{\beta}} = \underline{\smash{\gamma}}\,\underline{\smash{\beta}}^{'}$ and the subsums of the coordinates of $\underline{\alpha}^{'}$ and $\underline{\smash{\beta}}^{'}$ are pairwise distinct?
\end{problem}
\end{comment}

\section{The Main Lemma}
In this section we prove the main lemma which plays a crucial role in the proofs of our theorems.

\begin{lemma}\label{lemma_1}
Let $n\in \mathbb{Z}^{+}$. 
\begin{itemize}
\item[(i)] Let $C,D\in \mathbb{N}_{m}$, where $C = 
\{0,c_{2},c_{3}, \dots{}\}$ $(0 < c_{2} \le c_{3} \le \dots{})$, $|C|, |D| \ge 2$. If $\underline{\alpha}, \underline{\smash{\beta}} \in \mathbb{Z}^{+}_{m}$ satisfy $(\underline{\alpha}C) \cap [c_{2}n]_{m} = (\underline{\smash{\beta}}D) \cap [c_{2}n]_{m}$ and
$C \cap [c_{2}n] = D \cap [c_{2}n]$, \\
then $\underline{\alpha} \cap [n]_{m} = \underline{\smash{\beta}}\cap [n]_{m}$.
\item[(ii)] Let $\underline{\alpha}, \underline{\smash{\beta}} \in \mathbb{Z}^{+}_{m}$, where $\underline{\alpha} = (\alpha_{1}, \alpha_{2}, \dots{} ,\alpha_{s})$ and $C,D\in \mathbb{N}_{m}$ such that
$\underline{\alpha} \cap [\alpha_{1}n]_{m} = \underline{\smash{\beta}} \cap [\alpha_{1}n]_{m}$ and
$(\underline{\alpha}C) \cap [\alpha_{1}n] = (\underline{\smash{\beta}}D) \cap [\alpha_{1}n]$ hold. Then $C \cap [n]_{m} = D\cap [n]_{m}$.
\item[(iii)] Let $\underline{\alpha}, \underline{\alpha}^{'}, \underline{\smash{\beta}}, \underline{\smash{\beta}}^{'} \in \mathbb{Z}^{+}_{m}$, where $\underline{\alpha} = (\alpha_{1}, \alpha_{2}, \dots{} ,\alpha_{s})$. If $\underline{\alpha} \cap [\alpha_{1}n]_{m} = \underline{\smash{\beta}} \cap [\alpha_{1}n]_{m}$, \\ $\underline{\alpha}\,\underline{\alpha}^{'} \cap [\alpha_{1}n]_{m} = \underline{\smash{\beta}}\,\underline{\smash{\beta}}^{'} \cap [\alpha_{1}n]_{m}$. Then $\underline{\alpha}^{'}\cap [n]_{m} = \underline{\smash{\beta}}^{'}\cap [n]_{m}$.
\end{itemize}
\end{lemma}

\begin{proof}
(i): Assume that $\underline{\alpha}\cap [n]_{m} \ne \underline{\smash{\beta}}\cap [n]_{m}$. Let $k$ be the smallest integer with $\chi_{\underline{\alpha}}(k) \ne \chi_{\underline{\smash{\beta}}}(k)$. Since $C \cap [c_{2}n] = D \cap [c_{2}n]$, the smallest non zero element of $D$ is $c_{2}$ and $\chi_{C}(c_{2}) = \chi_{D}(c_{2})$. Consider the tuples $(j_{1}, j_{2}, \dots{} ,j_{s})$ with $kc_{2} = \sum_{i=1}^{s}\alpha_{i}c_{j_{i}}$. In this representation, if $\alpha_{i} > k$, then $c_{j_{i}} = 0$, that is $j_{i} = 1$ and similarly, in the representation $kc_{2} = \sum_{i=1}^{t}\beta_{i}d_{j_{i}}$, if 
$\beta_{i} > k$, then $d_{j_{i}} = 0$, that is $j_{i} = 1$. Let $K\in \mathbb{N}$ such that $\alpha_{K} < k$, $\beta_{K} < k$, but $\alpha_{K+1} \ge k$ or  $\alpha_{K+1}$ does not exist. Thus $\beta_{K+1} \ge k$ or  $\beta_{K+1}$ does not exist. If $kc_{2} = \sum_{i=1}^{s}\alpha_{i}c_{j_{i}}$ and there exists a $K < i$ with $c_{j_{i}} > 0$, then $\alpha_{i} = k$, $c_{j_{i}} = c_{2}$ and if $h\ne i$, then $c_{j_{h}} = 0$, that is $j_{h} = 1$.
Similarly, if $kc_{2} = \sum_{i=1}^{t}\beta_{i}d_{j_{i}}$ and there exists a $K < i$ with $d_{j_{i}} > 0$, then $\beta_{i} = k$, $d_{j_{i}} = c_{2}$ and if $h\ne i$, then $d_{j_{h}} = 0$, that is $j_{h} = 1$. It follows that 
\[
R_{C,\underline{\alpha}}(kc_{2}) = |\{(j_{1}, j_{2}, \dots{},j_{s}): kc_{2} = \sum_{i=1}^{s}\alpha_{i}c_{j_{i}}\}| 
\]
\[
= |\{(j_{1}, j_{2}, \dots{} ,j_{K}): kc_{2} = \sum_{i=1}^{K}\alpha_{i}c_{j_{i}}\}| +  \chi_{\underline{\alpha}}(k)\chi_{C}(c_{2}).
\]
Similarly, 
\[
R_{D,\underline{\smash{\beta}}}(kc_{2}) = |\{(j_{1}, j_{2}, \dots{} ,j_{K}): kc_{2} = \sum_{i=1}^{K}\beta_{i}d_{j_{i}}\}| +  \chi_{\underline{\smash{\beta}}}(k)\chi_{D}(c_{2}).
\]
Since $R_{C,\underline{\alpha}}(kc_{2}) = R_{D,\underline{\smash{\beta}}}(kc_{2})$, $C \cap [c_{2}n]_{m} = D \cap [c_{2}n]_{m}$ and $\underline{\alpha} \cap [k-1]_{m} = \underline{\smash{\beta}} \cap [k-1]_{m}$, then $kc_{2} = \sum_{i=1}^{K}\alpha_{i}c_{j_{i}}$ holds if and only if $kc_{2} = \sum_{i=1}^{K}\beta_{i}d_{j_{i}}$. It follows that
\[
|\{(j_{1}, j_{2}, \dots{} ,j_{K}): kc_{2} = \sum_{i=1}^{K}\alpha_{i}c_{j_{i}}\}| = |\{(j_{1}, j_{2}, \dots{} ,j_{K}): kc_{2} = \sum_{i=1}^{K}\beta_{i}d_{j_{i}}\}|,
\]
and so $\chi_{\underline{\alpha}}(k)\chi_{C}(c_{2}) = \chi_{\underline{\smash{\beta}}}(k)\chi_{D}(c_{2})$,
where $\chi_{C}(c_{2}) = \chi_{D}(c_{2}) > 0$, which implies that $\chi_{\underline{\alpha}}(k) = \chi_{\underline{\smash{\beta}}}(k)$, a contradiction.

(ii): Assume that $C\cap [n]_{m} \ne D\cap [n]_{m}$. Let $k$ be the smallest integer with $\chi_{C}(k) \ne \chi_{D}(k)$. Let $\underline{\smash{\beta}} = (\beta_{1}, \beta_{2}, \dots{} ,\beta_{t})$, where $\beta_{1} \le \beta_{2}\le \dots{} \le \beta_{t}$. Since $\underline{\alpha} \cap [\alpha_{1}n]_{m} = \underline{\smash{\beta}} \cap [\alpha_{1}n]_{m}$, thus $\alpha_{1} = \beta_{1}$ and $\chi_{\underline{\alpha}}(\alpha_{1}) = \chi_{\underline{\smash{\beta}}}(\alpha_{1})$. Consider the tuples $(j_{1}, j_{2}, \dots{} ,j_{s})$ with $\alpha_{1}k = \sum_{i=1}^{s}\alpha_{i}c_{j_{i}}$. 
%In this representation, if $\alpha_{i} > \alpha_{1}k$, then $c_{j_{i}} = 0$, that is $j_{i} = 1$ for every $i$. Similarly, in the representation $k\alpha_{1} = \sum_{i}\beta_{i}d_{j_{i}}$, if 
%$\beta_{i} > k\alpha_{1}$, then $d_{j_{i}} = 0$, that is $j_{i} = 1$ for every $i$.
Let $K\in \mathbb{N}$ such that $\alpha_{K} < k$, $\beta_{K} < k$, but $\alpha_{K+1} \ge k$ or  $\alpha_{K+1}$ does not exist. Thus $\beta_{K+1} \ge k$ or  $\beta_{K+1}$ does not exist.
%Since $\underline{\alpha} \cap [\alpha_{1}n]_{m} = \underline{\smash{\beta}} \cap [\alpha_{1}n]_{m}$,  
%if $\alpha_{1}k = \sum_{i=1}^{s}\alpha_{i}c_{j_{i}}$ is such that there exists an $i$ with $c_{j_{i}} \ge k$, then $\alpha_{i} = \alpha_{1}$, $c_{j_{i}} = k$ and if $h \ne i$ then $c_{j_{h}} = 0$ that is $j_{h} = 1$. 
%Similarly, if $\alpha_{1}k = \sum_{i=1}^{t}\beta_{i}d_{j_{i}}$ is such that there exists an $i$ with $d_{j_{i}} \ge k$, then $\beta_{i} = \alpha_{1}$, $d_{j_{i}} = k$ and if $h\ne i$, then $d_{j_{h}} = 0$, that is $j_{h} = 1$. It follows that 
Similar to the previous case we get
%\[
%R_{C,\underline{\alpha}}(\alpha_{1}k) = |\{(j_{1}, j_{2}, \dots{} j_{s}): \alpha_{1}k = \sum_{i=1}^{s}\alpha_{i}c_{j_{i}}\}| 
%\]
\[
R_{C,\underline{\alpha}}(\alpha_{1}k) = |\{(j_{1}, j_{2}, \dots{} ,j_{K}): \alpha_{1}k = \sum_{i=1}^{K}\alpha_{i}c_{j_{i}}, c_{j_{i}} < k \textnormal{ for every } i\}| +  \chi_{\underline{\alpha}}(\alpha_{1})\chi_{C}(k).
\]
and 
\[
R_{D,\underline{\smash{\beta}}}(\alpha_{1}k) = |\{(j_{1}, j_{2}, \dots{} ,j_{K}): \alpha_{1}k = \sum_{i=1}^{K}\beta_{i}d_{j_{i}}, d_{j_{i}} < k \textnormal{ for every } i\}| + \chi_{\underline{\beta}}(\alpha_{1})\chi_{D}(k).
\]
Since $R_{C,\underline{\alpha}}(\alpha_{1}k) = R_{D,\underline{\smash{\beta}}}(\alpha_{1}k)$, $\underline{\alpha} \cap [\alpha_{1}n]_{m} = \underline{\smash{\beta}} \cap [\alpha_{1}n]_{m}$ and $C \cap [k-1]_{m} = D \cap [k-1]_{m}$, then $\alpha_{1}k = \sum_{i=1}^{K}\alpha_{i}c_{j_{i}}$, where $c_{j_{i}} < k \textnormal{ for every } i$ holds if and only if $\alpha_{1}k = \sum_{i=1}^{K}\beta_{i}d_{j_{i}}$ with $d_{j_{i}} < k \textnormal{ for every } i$. Thus we have 
\[
|\{(j_{1}, j_{2}, \dots{} ,j_{K}): \alpha_{1}k = \sum_{i=1}^{K}\alpha_{i}c_{j_{i}}, c_{j_{i}} < k \textnormal{ for every } i\}| 
\]
\[
=
|\{(j_{1}, j_{2}, \dots{} ,j_{K}): \alpha_{1}k = \sum_{i=1}^{K}\beta_{i}d_{j_{i}}, d_{j_{i}} < k \textnormal{ for every } i\}|
\]
and so $\chi_{\underline{\alpha}}(\alpha_{1})\chi_{C}(k) = \chi_{\underline{\smash{\beta}}}(\alpha_{1})\chi_{D}(k)$,
where $\chi_{\underline{\alpha}}(\alpha_{1}) = \chi_{\underline{\smash{\beta}}}(\alpha_{1}) > 0$, which implies that $\chi_{C}(k) = \chi_{D}(k)$, which is a contradiction.

(iii): Assume that $\underline{\alpha}^{'}\cap [n]_{m} \ne \underline{\smash{\beta}}^{'}\cap [n]_{m}$. Let $k$ be the smallest integer with $\chi_{\underline{\alpha}^{'}}(k) \ne \chi_{\underline{\smash{\beta}}^{'}}(k)$.
Since $\underline{\alpha} \cap [\alpha_{1}n]_{m} = \underline{\smash{\beta}} \cap [\alpha_{1}n]_{m}$, then $\alpha_{1} = \beta_{1}$ and $\chi_{\underline{\alpha}}(\alpha_{1}) = \chi_{\underline{\smash{\beta}}}(\alpha_{1})$. It follows that in the equation  $\alpha_{1}k = \alpha_{i}\alpha^{'}_{j}$, $\alpha^{'}_{j} \ge k$ holds if 
$\alpha^{'}_{j} = k$, and $\alpha_{i} = \alpha_{1}$. Similarly, in $\alpha_{1}k = \beta_{i}\beta^{'}_{j}$, $\beta^{'}_{j} \ge k$ holds if 
$\beta^{'}_{j} = k$ and $\beta_{i} = \alpha_{1}$. Thus we have
\[
\chi_{\underline{\alpha}\,\underline{\alpha}^{'}}(\underline{\alpha}_{1}k) = 
|\{(i,j): \alpha_{1}k = \alpha_{i}\alpha^{'}_{j}, \alpha_{i}\in \underline{\alpha},\alpha_{j}^{'}\in \underline{\alpha}^{'}\}| 
\]
\[
= |\{(i,j): \alpha_{1}k = \alpha_{i}\alpha^{'}_{j}, \alpha^{'}_{j} < k, \alpha_{i}\in \underline{\alpha}, \alpha_{j}^{'}\in \underline{\alpha}^{'}\}| + \chi_{\underline{\alpha}}(\alpha_{1})\chi_{\underline{\alpha}^{'}}(k)
\]
and
\[
\chi_{\underline{\smash{\beta}}\,\underline{\smash{\beta}}^{'}}(\underline{\alpha}_{1}k) = |\{(i,j): \alpha_{1}k = \beta_{i}\beta^{'}_{j}, \beta^{'}_{j} < k,\beta_{i}\in \underline{\smash{\beta}},\beta_{j}^{'}\in \underline{\smash{\beta}}^{'} \}| + \chi_{\underline{\smash{\beta}}}(\alpha_{1})\chi_{\underline{\smash{\beta}}^{'}}(k).
\]
Since $\underline{\alpha} \cap [\alpha_{1}n]_{m} = \underline{\smash{\beta}} \cap [\alpha_{1}n]_{m}$,
$\underline{\alpha}^{'} \cap [k-1]_{m} = \underline{\smash{\beta}}^{'} \cap [k-1]_{m}$,
then $\alpha_{1}k = \alpha_{i}\alpha^{'}_{j}$, $\alpha_{i}\in \underline{\alpha},\alpha_{j}^{'}\in \underline{\alpha}^{'}$, $\alpha^{'}_{j} < k$ holds if and only if
 $\alpha_{1}k = \beta_{i}\beta^{'}_{j}$, $\beta_{i}\in \underline{\smash{\beta}}$, $\beta_{j}^{'}\in \underline{\smash{\beta}}^{'}$, $\beta^{'}_{j} < k$ and so
 \[
 |\{(i,j): \alpha_{1}k = \alpha_{i}\alpha^{'}_{j}, \alpha_{i}\in \underline{\alpha}, \alpha_{j}^{'}\in \underline{\alpha}^{'}, \alpha^{'}_{j} < k\}| =  |\{(i,j): \alpha_{1}k = \beta_{i}\beta^{'}_{j}, \beta_{i}\in \underline{\smash{\beta}},\beta_{j}^{'}\in \underline{\smash{\beta}}^{'}, \beta^{'}_{j} < k\}|.
 \]
Since $\chi_{\underline{\alpha}\,\underline{\alpha}^{'}}(\underline{\alpha}_{1}k) = \chi_{\underline{\smash{\beta}}\,\underline{\smash{\beta}}^{'}}(\underline{\alpha}_{1}k)$, thus
$\chi_{\underline{\alpha}}(\alpha_{1})\chi_{\underline{\alpha}^{'}}(k) = \chi_{\underline{\smash{\beta}}}(\alpha_{1})\chi_{\underline{\smash{\beta}}^{'}}(k)$
and $\chi_{\underline{\alpha}}(\alpha_{1}) = \chi_{\underline{\smash{\beta}}}(\alpha_{1}) > 0$, implies that $\chi_{\underline{\alpha}^{'}}(k) = \chi_{\underline{\smash{\beta}}^{'}}(k)$, a contradiction.
\end{proof}

\section{Proof of Theorem \ref{prop_1}}

(1): Since $R_{A,\underline{\alpha}}(n) = R_{B,\underline{\smash{\beta}}}(n)$ for every $n$, the smallest value of $n$ for which $R_{A,\underline{\alpha}}(n) > 0$ is $n = \sum_{i=1}^{s}\alpha_{i}a_{1} = \sum_{i=1}^{t}\beta_{i}b_{1}$.
Since
$R_{A+u,\underline{\alpha}}(n) = R_{B+v,\underline{\smash{\beta}}}(n)$ for every $n$, and 
the smallest value of $n$ for which $R_{A+u,\underline{\alpha}}(n) > 0$ is $n = \sum_{i=1}^{s}(a_{1}+u)\alpha_{i} = \sum_{i=1}^{t}(b_{1}+v)\beta_{i}$. It follows that 
\[
\sum_{i=1}^{s}(a_{1}+u)\alpha_{i} = a_{1}\sum_{i=1}^{s}\alpha_{i} + u\sum_{i=1}^{s}\alpha_{i} = \sum_{i=1}^{t}(b_{1}+v)\beta_{i} = b_{1}\sum_{i=1}^{t}\beta_{i} + v\sum_{i=1}^{t}\beta_{i}
\]
and so
\[
u\sum_{i=1}^{s}\alpha_{i} = v\sum_{i=1}^{t}\beta_{i}.
\]
Let $d$ be the greatest common divisor of $\sum_{i=1}^{s}\alpha_{i}$ and $\sum_{i=1}^{t}\beta_{i}$. Since
\[
u\frac{\sum_{i=1}^{s}\alpha_{i}}{d} = v\frac{\sum_{i=1}^{t}\beta_{i}}{d},
\]
it follows that $\frac{\sum_{i=1}^{t}\beta_{i}}{d}$ divides $u$, so there exists an integer $w$ such that

\begin{equation}\label{eq:5}
u = w\frac{\sum_{i=1}^{t}\beta_{j}}{\textnormal{ gcd}(\sum_{i=1}^{s}\alpha_{i}, \sum_{i=1}^{t}\beta_{i})}
\textnormal{ and }
v = w\frac{\sum_{i=1}^{s}\alpha_{i}}{\textnormal{ gcd}(\sum_{i=1}^{s}\alpha_{i}, \sum_{i=1}^{t}\beta_{i})}.
\end{equation}

Conversely, if there is an integer $w$ such that \eqref{eq:5} holds, then 
\[
R_{A+u,\underline{\alpha}}(n) =  R_{A,\underline{\alpha}}\left(n-\sum_{i=1}^{s}\alpha_{i}u\right) =  R_{A,\underline{\alpha}}\left(n-w\frac{\sum_{i=1}^{s}\alpha_{i}\sum_{i=1}^{t}\beta_{i}}{d}\right) 
\]
\[
= R_{B,\underline{\smash{\beta}}}\left(n-w\frac{\sum_{i=1}^{s}\alpha_{i}\sum_{i=1}^{t}\beta_{i}}{d}\right) = R_{B,\underline{\smash{\beta}}}\left(n-\sum_{i=1}^{t}\beta_{i}v\right) = R_{B+v,\underline{\smash{\beta}}}(n).
\]

(2): The second statement follows from the first part in the following way. 
Since $R_{A,\underline{\alpha}}(n) = R_{B,\underline{\smash{\beta}}}(n)$ for every $n$, the smallest value of $n$ for which $R_{A,\underline{\alpha}}(n) > 0$ is $n = \sum_{i=1}^{s}\alpha_{i}a_{1} = \sum_{i=1}^{t}\beta_{i}b_{1}$. It follows that
\[
a_{1}\frac{\sum_{i=1}^{s}\alpha_{i}}{d} = b_{1}\frac{\sum_{i=1}^{t}\beta_{i}}{d},
\]
where $d = \textnormal{ gcd}(\sum_{i=1}^{s}\alpha_{i}, \sum_{i=1}^{t}\beta_{i})$ and so there exists a nonnegative integer $w$ such that
\[
a_{1} = w\frac{\sum_{i=1}^{t}\beta_{i}}{d}, \hspace*{10mm} b_{1} 
= w\frac{\sum_{i=1}^{s}\alpha_{i}}{d}.
\]
It follows from the previous part that $R_{A-a_{1},\underline{\alpha}} = R_{B-b_{1},\underline{\smash{\beta}}}$.

\section{Proof of Theorem \ref{theorem_1}}
(i): Let $n\in \mathbb{Z}^{+}$. Since $\underline{\alpha} A \cap [\underline{\alpha}_{1}n]_{m} = \underline{\alpha}B \cap [\underline{\alpha}_{1}n]_{m}$ and $\underline{\alpha} \cap [\underline{\alpha}_{1}n]_{m} = \underline{\alpha} \cap [\underline{\alpha}_{1}n]_{m}$, then by (ii) of Lemma \ref{lemma_1}, $A \cap [n]_{m} = B \cap [n]_{m}$ which implies that $A = B$. \\
(ii): Let $n\in \mathbb{Z}^{+}$ and $A = \{0, a_{2}, \dots{}\}$, $0 < a_{2} \le \dots{}$. Since $\underline{\alpha} A \cap [a_{2}n]_{m} = \underline{\smash{\beta}}A \cap [a_{2}n]_{m}$, $A \cap [a_{2}n]_{m} = A \cap [a_{2}n]_{m}$, then by (i) of Lemma \ref{lemma_1}, $ \underline{\alpha}\cap [n]_{m} = \underline{\smash{\beta}} \cap [n]_{m}$ which implies that $\underline{\alpha} = \underline{\smash{\beta}}$.

\subsection{Proof of Corollary 1.7}
(i): Obviously, we can assume that $0 \notin A\cap B$. If 
$R_{A,\underline{\alpha}} = R_{B,\underline{\smash{\beta}}}$, then by (2) in Proposition \ref{prop_1}, $R_{A-a_{1},\underline{\alpha}} = R_{B-b_{1},\underline{\alpha}}$. It follows from (i) in Theorem \ref{theorem_1} that
$A-a_{1} = B-b_{1}$. It is clear that $\sum_{i=1}^{s}\alpha_{i}a_{1}$ is the smallest positive integer for which 
$R_{A,\underline{\alpha}}$ is positive and $\sum_{i=1}^{s}\alpha_{i}b_{1}$ is the smallest positive integer such that 
$R_{B,\underline{\alpha}}$ is positive, then $\sum_{i=1}^{s}\alpha_{i}a_{1} = \sum_{i=1}^{s}\alpha_{i}b_{1}$. Thus we have $a_{1} = b_{1}$, and then $A-a_{1} = B-b_{1} = B-a_{1}$ which implies that $A = B$.\\ 
(ii): If $R_{A,\underline{\alpha}} = R_{A,\underline{\smash{\beta}}}$, then by (2) in Proposition \ref{prop_1}, $R_{A-a_{1},\underline{\alpha}} = R_{A-a_{1},\underline{\smash{\beta}}}$, where $0 \in A-a_{1}$. It follows from (ii) of Theorem \ref{theorem_1} that $\underline{\alpha} = \underline{\smash{\beta}}$.

\section{Proof of Theorem \ref{theorem_2}}

The proof is divided into three steps. In the first step for $i = 1, 2, \dots{}$ we define the vectors and the multisets

\begin{equation}\label{eq:6}
\underline{\alpha}^{(i)}_{1}, \underline{\alpha}^{(i)}_{2}, \underline{\smash{\beta}}^{(i)}_{1},  
\underline{\smash{\beta}}^{(i)}_{2}, \underline{\smash{\gamma}}^{(i)}, \underline{\delta}^{(i)}\in \mathbb{Z}_{m}^{+}, 
\hspace*{3mm}
C^{(i)}, D^{(i)}\in \mathbb{N}
\end{equation}
which satisfy the following equations.

\begin{equation}\label{eq:7}
A = \underline{\alpha}^{(1)}_{1}\underline{\alpha}^{(2)}_{1}\dots{} \underline{\alpha}^{(i)}_{1}C^{(i)}, 
\end{equation}

\begin{equation}\label{eq:8}
B = \underline{\smash{\beta}}^{(1)}_{1} \underline{\smash{\beta}}^{(2)}_{1} \dots{}  \underline{\smash{\beta}}^{(i)}_{1}D^{(i)},
\end{equation}

\begin{equation}\label{eq:9}
\underline{\alpha} = \underline{\alpha}^{(1)}_{2}\underline{\alpha}^{(2)}_{2} \dots{}  \underline{\alpha}^{(i)}_{2}\underline{\smash{\gamma}}^{(i)}, 
\end{equation}

\begin{equation}\label{eq:10}
\underline{\smash{\beta}} =  \underline{\smash{\beta}}^{(1)}_{2}\underline{\smash{\beta}}^{(2)}_{2} \dots{}  \underline{\smash{\beta}}^{(i)}_{2}\underline{\delta}^{(i)},
\end{equation}
\[
\underline{\alpha}^{(i)}_{1}\underline{\alpha}^{(i)}_{2} = \underline{\smash{\beta}}^{(i)}_{1} \underline{\smash{\beta}}^{(i)}_{2} := (u^{(i)}_{1}, u^{(i)}_{2}, \dots{} ).
\]
For $i = 1, 2, \dots{}$ we choose vectors $\underline{\alpha}^{(i)}_{1}, \underline{\alpha}^{(i)}_{2}, \underline{\smash{\beta}}^{(i)}_{1}, \underline{\smash{\beta}}^{(i)}_{2}$ with $u^{(i)}_{1} \ge 2$. We show that bounded number of vectors can be chosen with this property. Let us suppose that the maximal number of $4$-tuples ($\underline{\alpha}^{(i)}_{1}, \underline{\alpha}^{(i)}_{2}, \underline{\smash{\beta}}^{(i)}_{1}, \underline{\smash{\beta}}^{(i)}_{2}$) is $i_{0}$.
Now we fix $i_{0}$ and the corresponding vectors. Next, we choose 
$\underline{\alpha}^{(i)}_{1}, \underline{\alpha}^{(i)}_{2}, \underline{\smash{\beta}}^{(i)}_{1}, \underline{\smash{\beta}}^{(i)}_{2}$
with $(u^{(i)}_{1}, u^{(i)}_{2}, \dots{} ) = (1,1,\dots{})$ for $i = i_{0}+1, i_{0}+2, \dots{}$. We show that there exist bounded number of such vectors. Let us denote by $i_{1}$ the maximum number of these vectors.
Now we fix $i_{1}$ and the corresponding vectors. Let us suppose that we already defined $i_{0}, \dots{} i_{n-1}$. We show that for $i = i_{0}+i_{1}+\dots{} + i_{n-1}+1, i_{0}+ i_{1}+\dots{} + i_{n-1}+2, \dots{}$ there exist bounded number of such vectors $\underline{\alpha}^{(i)}_{1}, \underline{\alpha}^{(i)}_{2}, \underline{\smash{\beta}}^{(i)}_{1}, \underline{\smash{\beta}}^{(i)}_{2}$ with $(u^{(i)}_{1}, u^{(i)}_{2}, \dots{} ) = (1,n,\dots{})$. 
Let us suppose that the maximal number of $4$-tuples $(\underline{\alpha}^{(i)}_{1}, \underline{\alpha}^{(i)}_{2}, \underline{\smash{\beta}}^{(i)}_{1}, \underline{\smash{\beta}}^{(i)}_{2})$ is $i_{n}$. 
In the $n$-th step we fix $i_{n}$ and the corresponding $i_{n}$ vectors.

%We define the limits of the above  multisets in the following way. For a sequence of %multisets $C^{(1)}, C^{(2)}, \dots{}$, we write
%\[
%\lim_{n\rightarrow \infty}C^{(n)} = C
%\]
%if for every positive integer $n$, there exists a positive integer $N$ such that for %every $m \ge N$, then $C^{(m)} \cap [n]_{m} = C\cap [n]_{m}$.

In the second step, we prove that the following limits exist.

\begin{equation}\label{eq:11}
\lim_{n\rightarrow \infty}\underline{\alpha}^{(1)}_{1}\underline{\alpha}^{(2)}_{1}\dots{} \underline{\alpha}^{(i_{0}+i_{1}+\dots{}+i_{n})}_{1} = \underline{\alpha}_{1}   
\end{equation}

\begin{equation}\label{eq:12}
\lim_{n\rightarrow \infty}\underline{\smash{\beta}}^{(1)}_{1}\underline{\smash{\beta}}^{(2)}_{1}\dots{} \underline{\smash{\beta}}^{(i_{0}+i_{1}+\dots{}+i_{n})}_{1} = \underline{\smash{\beta}}_{1}    
\end{equation}

\begin{equation}\label{eq:13}
\lim_{n\rightarrow \infty}\underline{\alpha}^{(1)}_{2}\underline{\alpha}^{(2)}_{2}\dots{} \underline{\alpha}^{(i_{0}+i_{1}+\dots{}+i_{n})}_{2} = \underline{\alpha}_{2}   
\end{equation}

\begin{equation}\label{eq:14}
\lim_{n\rightarrow \infty}\underline{\smash{\beta}}^{(1)}_{2}\underline{\smash{\beta}}^{(2)}_{2}\dots{} \underline{\smash{\beta}}^{(i_{0}+i_{1}+\dots{}+i_{n})}_{2} = \underline{\smash{\beta}}_{2}    
\end{equation}

\begin{equation}\label{eq:15}
\lim_{n\rightarrow \infty} \underline{\smash{\gamma}}^{(i_{0}+i_{1}+\dots{}+i_{n})} = \underline{\smash{\gamma}}   
\end{equation}

\begin{equation}\label{eq:16}
\lim_{n\rightarrow \infty}\underline{\delta}^{(i_{0}+i_{1}+\dots{}+i_{n})} = \underline{\delta}   
\end{equation}

\begin{equation}\label{eq:17}
\lim_{n\rightarrow \infty}C^{(i_{0}+i_{1}+\dots{}+i_{n})} = C
\end{equation}

\begin{equation}\label{eq:18}
\lim_{n\rightarrow \infty}D^{(i_{0}+i_{1}+\dots{}+i_{n})} = D,  
\end{equation}
where $A = \underline{\smash{\alpha_{1}}}C$, 
$\underline{\alpha} = \underline{\smash{\gamma}}\,\underline{\smash{\alpha_{2}}}$, $B = \underline{\smash{\beta}}_{1}D$,
$\underline{\smash{\beta}} = \underline{\delta}\,\underline{\smash{\beta}}_{2}$ and $\underline{\alpha}_{1}\,\underline{\alpha}_{2} = \underline{\smash{\beta}}_{1} \underline{\smash{\beta}}_{2}$, thus we have $\underline{\alpha}A = (\underline{\smash{\alpha}}_{1}\,\underline{\smash{\alpha}}_{2})(\underline{\smash{\gamma}}C), 
\underline{\smash{\beta}}B = (\underline{\smash{\beta}}_{1}\,\underline{\smash{\beta}}_{2})(\underline{\delta}D)$. It follows from Theorem \ref{theorem_1} that $\underline{\smash{\gamma}}C = \underline{\delta}D$.

In the last step, we prove that $\underline{\smash{\gamma}}C = \underline{\delta}D$ is primitive.

We denote the coordinates of the vectors and components of the set of multisets defined in \eqref{eq:6}-\eqref{eq:10} in the following way. If $\underline{\smash{\mu}}^{(i)}_{j} \in \{\underline{\alpha}^{(i)}_{1}, \underline{\alpha}^{(i)}_{2}, \underline{\smash{\beta}}^{(i)}_{1}, \underline{\smash{\beta}}^{(i)}_{2}\}$, then $\underline{\smash{\mu}}^{(i)}_{j} = (\mu^{(i)}_{j,1}, \mu^{(i)}_{j,2}, \dots{})$. Moreover, if $\underline{\nu}^{(i)}\in \{\underline{\smash{\gamma}}^{(i)}, \underline{\delta}^{(i)}\}$, then 
$\underline{\nu}^{(i)} = (\nu^{(i)}_{1}, \nu^{(i)}_{2}, \dots{})$.
Furthermore, if $S^{(i)}\in \{C^{(i)}, D^{(i)}\}$, then $S^{(i)} = \{s^{(i)}_{1}, s^{(i)}_{2}, \dots{}\}$, where $s^{(i)}_{1} \le s^{(i)}_{2} \le \dots{}$.

Now we prove the first step.
Let $N$ be the smallest positive integer with $R_{A,\underline{\alpha}}(N) > 0$.
Assume that there exist vectors
$\underline{\alpha}^{'}, \underline{\alpha}^{''}, \underline{\smash{\beta}}^{'}$,  $\underline{\smash{\beta}}^{''}, \underline{\smash{\gamma}}^{''}, \underline{\delta}^{''}$ and multisets $C^{'}, D^{'}$ such that $\underline{\alpha}^{'}  \underline{\alpha}^{''} = \underline{\smash{\beta}}^{'}\,\underline{\smash{\beta}}^{''} := (u_{1}, u_{2}, \dots{} )$,  $1 < u_{1}$, $A = \underline{\alpha}^{'}C^{'}$, $B = \underline{\smash{\beta}}^{'}D^{'}$, $\underline{\alpha} = \underline{\alpha}^{''}\,\underline{\smash{\gamma}}^{''}, \underline{\smash{\beta}} = \underline{\smash{\beta}}^{''} \underline{\delta}^{''}$.

Then if there exist vectors and multisets defined in \eqref{eq:6} for some $1 \le i \le s$ which satisfy \eqref{eq:7}-\eqref{eq:10} then
\[
\underline{\alpha} A = \underline{\alpha}^{(1)}_{1}  \underline{\alpha}^{(2)}_{1} \underline{\alpha}^{(1)}_{2} \underline{\alpha}^{(2)}_{2}  \dots{}  \underline{\alpha}^{(1)}_{s} \underline{\alpha}^{(2)}_{s}\underline{\smash{\gamma}}^{(s)}C^{(s)}
\]
and we have $N =  u_{1,1} u_{2,1} \cdots{} u_{s,1} \gamma^{(s)}_{1} c^{(s)}_{2} \ge 2^{s}$, which implies that $s \le \log_{2}N$. It follows that there exists a positive integer $i_{0}$ such that for $1 \le j \le i_{0}$, there exist vectors and multisets defined in \eqref{eq:6}, but for $i_{0} + 1$, there do not exist such vectors and multisets. If there does not exist vectors
$\underline{\alpha}^{'}, \underline{\alpha}^{''}, \underline{\smash{\beta}}^{'}$,  $\underline{\smash{\beta}}^{''}, \underline{\smash{\gamma}}^{''}, \underline{\delta}^{''}$ and multisets $C^{'}, D^{'}$ such that $\underline{\alpha}^{'}  \underline{\alpha}^{''} = \underline{\smash{\beta}}^{'} \underline{\smash{\beta}}^{''} := (u_{1}, u_{2}, \dots{} ,u_{n})$, $1 < u_{1}$, $A = \underline{\alpha}^{'}C^{'}$, $B = \underline{\smash{\beta}}^{'}D^{'}$, $\underline{\alpha} = \underline{\alpha}^{''} \underline{\smash{\gamma}}^{''}, \underline{\smash{\beta}} = \underline{\smash{\beta}}^{''} \underline{\delta}^{''}$, then let $i_{0} = 0$ and $\underline{\alpha}^{(1)}_{1}  \underline{\alpha}^{(1)}_{2} \dots{} \underline{\alpha}^{(1)}_{i_{0}} = (1)$, 
$\underline{\alpha}^{(2)}_{1}  \underline{\alpha}^{(2)}_{2} \dots{} \underline{\alpha}^{(2)}_{i_{0}} = (1)$, $\underline{\alpha}^{(1)}_{1}  \underline{\alpha}^{(2)}_{1} \underline{\alpha}^{(1)}_{2} \underline{\alpha}^{(2)}_{2}  \dots{}  \underline{\alpha}^{(1)}_{i_{0}}\underline{\alpha}^{(2)}_{i_{0}} = (1)$. 
We fix these tuples and vectors defined above. Let
\[
\underline{\alpha}^{(1)}_{1} \underline{\alpha}^{(2)}_{1} \underline{\alpha}^{(1)}_{2} \underline{\alpha}^{(2)}_{2}  \dots{}  \underline{\alpha}^{(1)}_{i_{0}} \underline{\alpha}^{(2)}_{i_{0}} = (U_{1},U_{2}, \dots{}).
\]

Assume that if $i = i_{0} + j$, $1 \le j \le p$, then one can find vectors and multisets defined in \eqref{eq:6}-\eqref{eq:10} which satisfy $\underline{\alpha}^{(1)}_{i} \underline{\alpha}^{(2)}_{i} = (1,1,\dots{})$. 
Then we have
\[
\underline{\alpha}^{(1)}_{1}  \underline{\alpha}^{(2)}_{1}\underline{\alpha}^{(1)}_{2}\underline{\alpha}^{(2)}_{2}   \dots{}  \underline{\alpha}^{(1)}_{i_{0}} \underline{\alpha}^{(2)}_{i_{0}} \underline{\alpha}^{(1)}_{i_{0}+1} \underline{\alpha}^{(2)}_{i_{0}+1} \underline{\alpha}^{(1)}_{i_{0}+2} \underline{\alpha}^{(2)}_{i_{0}+2}   \dots{}  \underline{\alpha}^{(1)}_{i_{0}+p} \underline{\alpha}^{(2)}_{i_{0}+p}
\]
\[
= (U_{1}, U_{2}, \dots{}) \underbrace{(1, 1, \dots{}) \dots{}  (1, 1, \dots{})}_{p} = (\underbrace{U_{1}, U_{1}, \dots{} ,U_{1}}_{2^{p}}, \dots{}).
\]
It follows from 
\[
\underline{\alpha} A = \underline{\alpha}^{(1)}_{1}  \underline{\alpha}^{(2)}_{1} \underline{\alpha}^{(1)}_{2} \underline{\alpha}^{(2)}_{2}  \dots{}  \underline{\alpha}^{(1)}_{i_{0}+p} \underline{\alpha}^{(2)}_{i_{0}+p}\underline{\smash{\gamma}}^{(i_{0}+p)}C^{(i_{0}+p)},
\]
that
\[
N = U_{1}\cdot \gamma^{(i_{0}+p)}_{1}\cdot c^{(i_{0}+p)}_{2},
\]
and
\[
R_{A,\underline{\alpha}}(N) \ge \chi_{C^{(i_{0}+p)}}(c^{(i_{0}+p)}_{2})\cdot \chi_{\underline{\smash{\gamma}}^{(i_{0}+p)}}(\smash{\gamma}^{(i_{0}+p)}_{1})\cdot \chi_{\underline{\alpha}^{(1)}_{1}  \underline{\alpha}^{(2)}_{1} \underline{\alpha}^{(1)}_{2} \underline{\alpha}^{(2)}_{2}  \dots{} \underline{\alpha}^{(1)}_{i_{0}+p} \underline{\alpha}^{(2)}_{i_{0}+p}}(U_{1}) \ge 2^{p},
\]
and so $p \le \log_{2}R_{A,\underline{\alpha}}(N)$. It follows that there exists a natural number $i_{1}$ such that for $i_{0} + 1 \le i \le i_{0} + i_{1}$ one can find vectors and multisets of the form \eqref{eq:6}, but for $i_{1} + 1$ there do not exist such vectors and multisets. We fix these tuples and vectors.

Now suppose that we already have the natural numbers $i_{1}, \dots{} ,i_{n-1}$ and the corresponding vectors and multisets defined in \eqref{eq:6}. Assume that for $i = i_{0} + i_{1} + \dots{} + i_{n-1} + j$, $1 \le j \le r$, one could find $r$ vectors and multisets defined in \eqref{eq:6} with properties \eqref{eq:7}-\eqref{eq:10} and 

\begin{equation}\label{eq:19}
    \underline{\alpha}^{(1)}_{i_{0} + i_{1} + \dots{} + i_{n-1} + j}  \underline{\alpha}^{(2)}_{i_{0} + i_{1} + \dots{} + i_{n-1} + j} = \underline{\smash{\beta}}^{(1)}_{i_{0} + i_{1} + \dots{} + i_{n-1} + j}  \underline{\smash{\beta}}^{(2)}_{i_{0} + i_{1} + \dots{} + i_{n-1} + j} = (1,n,\dots{})
\end{equation}
Let $b = i_{0} + i_{1} + \dots{} + i_{n-1} + r$.
Then, $\underline{\alpha}A = \underline{\alpha}^{(1)}_{1}\underline{\alpha}^{(1)}_{2}\underline{\alpha}^{(2)}_{1}\underline{\alpha}^{(2)}_{2} \cdots{} \underline{\alpha}^{(b)}_{1}\underline{\alpha}^{(b)}_{2}\underline{\smash{\gamma}}^{(b)}C^{(b)}$, $N = c^{(b)}_{2} \cdot \gamma^{(b)}_{1} \cdot U_{1}$ and
$Nn = c^{(b)}_{2} \cdot \gamma^{(b)}_{1} \cdot U_{1} \cdot n$,
\[
R_{A,\underline{\alpha}}(Nn) 
= R_{C^{(b)}, \underline{\alpha}^{(1)}_{1}  \underline{\alpha}^{(2)}_{1} \underline{\alpha}^{(1)}_{2} \underline{\alpha}^{(2)}_{2}   \dots{}    \underline{\alpha}^{(1)}_{b}\underline{\alpha}^{(2)}_{b}
\underline{\smash{\gamma}}^{(b)}}
(c^{(b)}_{2} \cdot \gamma^{(b)}_{1} \cdot U_{1} \cdot n)
\]
\[
\ge \chi_{C^{(b)}}(c^{(b)}_{2})\cdot \chi_{\underline{\smash{\gamma}}^{(b)}}(\gamma^{(b)}_{1})\cdot \chi_{\underline{\alpha}^{(1)}_{1} \underline{\alpha}^{(2)}_{1}\underline{\alpha}^{(1)}_{2} \underline{\alpha}^{(2)}_{2}  \dots{} \underline{\alpha}^{(1)}_{b} \underline{\alpha}^{(2)}_{b}}(U_{1}\cdot n).
\]
On the other hand,
\[
\underline{\alpha}^{(1)}_{1}  \underline{\alpha}^{(1)}_{2} \underline{\alpha}^{(2)}_{1} \underline{\alpha}^{(2)}_{2} \dots{}   
\dots{} \underline{\alpha}^{(i_{0})}_{1}  \underline{\alpha}^{(i_{0})}_{2} \underline{\alpha}^{(i_{0}+1)}_{1} \underline{\alpha}^{(i_{0}+1)}_{2} \dots{} 
\underline{\alpha}^{(i_{0} + i_{1} + \dots{} + i_{n-1})}_{1} \underline{\alpha}^{(i_{0} + i_{1} + \dots{} + i_{n-1})}_{2}
\]
\[
\underline{\alpha}^{(i_{0} + i_{1} + \dots{} + i_{n-1} + 1)}_{1} \underline{\alpha}^{(i_{0} + i_{1} + \dots{} + i_{n-1} + 1)}_{2}
\cdots{}
\underline{\alpha}^{(b)}_{1} \underline{\alpha}^{(b)}_{2}
\]
\[
= (U_{1}, U_{2}, \dots{})(1,  \dots{}) \dots{} (1, \dots{})\underbrace{(1,n, \dots{} )\cdots{}(1,n, \dots{} )}_{r} 
\]
\[
= (U_{1}, U_{2}, \dots{})(1, \dots{}) \dots{} (1, \dots{})(1, \underbrace{n, \dots{} ,n}_{r}, \dots{})
\]
\[
= (\dots{}, \underbrace{nU_{1}, nU_{1}, \dots{} ,nU_{1}}_{r}, \dots{} )
\]
so 
\[
\chi_{\underline{\alpha}^{(1)}_{1}  \underline{\alpha}^{(2)}_{1} \underline{\alpha}^{(1)}_{2} \underline{\alpha}^{(2)}_{2}   \dots{} \underline{\alpha}^{(1)}_{b} \underline{\alpha}^{(2)}_{b}}(U_{1}\cdot n) \ge r.
\]
This implies that $R_{A,\underline{\alpha}}(Nn) \ge r$ and so there exists a natural number $i_{n}$ and for with $i_{0} + i_{1} + \dots{} + i_{n-1} + 1 \le i \le i_{0} + i_{1} + \dots{} + i_{n-1} + i_{n}$ there exists vectors and multisets defined in \eqref{eq:6}
with properties \eqref{eq:7}-\eqref{eq:10} and \eqref{eq:19} such that $i_{n} + 1$ such vectors and multisets already do not exist. We fix these tuples  and multisets.

In the second step, we can assume that $i_{0} + i_{1} + \dots{} + i_{n}$ tends to infinity as 
$n \rightarrow \infty$. If $l > i_{0} + i_{1} + \dots{} + i_{n}$, then 
$\underline{\alpha}^{(l)}_{1} \underline{\alpha}^{(l)}_{2} = \underline{\smash{\beta}}^{(l)}_{1} \underline{\smash{\beta}}^{(l)}_{2} = (1, u^{(l)}_{2}, \dots{})$,
where $u^{(l)}_{2} > n$, and so $\underline{\alpha}^{(l)}_{1,1} =  \underline{\alpha}^{(l)}_{2,1} = \underline{\smash{\beta}}^{(l)}_{1,1} = \underline{\smash{\beta}}^{(l)}_{2,1} = 1$,  $\alpha^{(l)}_{1,2} > n$, $\alpha^{(l)}_{2,2} > n$, $\beta^{(l)}_{1,2} > n$, $\beta^{(l)}_{2,2} > n$. It follows that $\underline{\alpha}^{(1)}_{1} \underline{\alpha}^{(1)}_{2} \dots{} \underline{\alpha}^{(l)}_{1} \cap [n]_{m} = 
\underline{\alpha}^{(1)}_{1} \underline{\alpha}^{(2)}_{1} \dots{} \underline{\alpha}^{(i_{0} + i_{1} + \dots{} + i_{n})}_{1}\cap [n]_{m}$, which implies that the limit \eqref{eq:11} exists.
One can see that the limits \eqref{eq:12}, \eqref{eq:13} and \eqref{eq:14} exist in the same way. 
Let $\underline{\alpha} = (\alpha_{1}, \alpha_{2}, \dots{})$. 
If $l > i_{0} + i_{1} + \dots{} + i_{\alpha_{1}n}$, then $\underline{\alpha}^{(1)}_{1}  \dots{} \underline{\alpha}^{(l)}_{1}\cap [\alpha_{1}n]_{m} = \underline{\alpha}^{(1)}_{1}  \dots{} \underline{\alpha}^{(i_{0} + i_{1} + \dots{} + i_{\alpha_{1}n})}_{1}\cap [\alpha_{1}n]_{m}$.
It follows from \eqref{eq:9} and (iii) of the Lemma \ref{lemma_1} that $\underline{\gamma}^{(l)}\cap [n]_{m} = \underline{\gamma}^{(i_{0} + i_{1} + \dots{} + i_{\alpha_{1}n})}\cap [n]_{m}$, so the limit \eqref{eq:15} exists. The limit \eqref{eq:16} can be proved in a similar way.

If $l > i_{0} + i_{1} + \dots{} + i_{\alpha_{1}n}$, then $\underline{\alpha}^{(1)}_{1}  \dots{} \underline{\alpha}^{(l)}_{1}  \cap [\alpha_{1}n]_{m} = 
\underline{\alpha}^{(1)}_{1}  \dots{} \underline{\alpha}^{(i_{0} + i_{1} + \dots{} + i_{\alpha_{1}n})} \cap [\alpha_{1}n]_{m}$ and 
\[
A = \underline{\alpha}^{(1)}_{1}  \dots{} \underline{\alpha}^{(l)}_{1}C^{(l)} 
= \underline{\alpha}^{(1)}_{1}  \dots{} \underline{\alpha}^{(i_{0} + i_{1} + \dots{} + i_{\alpha_{1}n})}_{1}C^{(i_{0} + i_{1} + \dots{} + i_{\alpha_{1}n})} 
\]
implies that
\[
\underline{\alpha}^{(1)}_{1}  \dots{} \underline{\alpha}^{(l)}_{1}C^{(l)}  \cap [\alpha_{1}n]_{m} = 
\underline{\alpha}^{(1)}_{1}  \dots{} \underline{\alpha}^{(i_{0} + i_{1} + \dots{} + i_{\alpha_{1}n})}_{1}C^{(i_{0} + i_{1} + \dots{} + i_{\alpha_{1}n})} \cap [\alpha_{1}n]_{m}.
\]
From (ii) in Lemma \ref{lemma_1}, we get that $C^{(l)} \cap [n]_{m} = C^{(i_{0} + i_{1} + \dots{} + i_{\alpha_{1}n})} \cap [n]_{m}$, thus the limit \eqref{eq:17} exists. It can be proved similarly that the limit \eqref{eq:18} is also exists.

Now we show that $\underline{\smash{\gamma}}C = \underline{\delta}D$ holds.
Let $n$ be a positive integer.
If $l$ is large enough, then $C^{(l)}\cap [n]_{m} = C\cap [n]_{m}$, $\underline{\alpha}_{1} \cap [n]_{m} = \underline{\alpha}^{(1)}_{1}  \dots{} \underline{\alpha}^{(l)}_{1}\cap [n]_{m}$ and so
\[
A \cap [n]_{m} = \underline{\alpha}^{(1)}_{1}  \dots{} \underline{\alpha}^{(l)}_{1}C^{(l)}  \cap [n]_{m} = ((\underline{\alpha}^{(1)}_{1}  \dots{} \underline{\alpha}^{(l)}_{1}\cap [n]_{m})(C^{(l)} \cap [n]_{m}))\cap [n]_{m} 
\]
\[
= ((\underline{\alpha}_{1}\cap [n]_{m})(C\cap [n]_{m}))\cap [n]_{m} = (\underline{\alpha}_{1}C)\cap[n]_{m},
\]
which implies that $A = \underline{\alpha}_{1}C$. One can prove similarly that $\underline{\alpha} = \underline{\smash{\gamma}}\,\underline{\alpha}_{2}$,
$\underline{\smash{\beta}} = \underline{\delta}\,\underline{\smash{\beta}}_{2}$, $B = \underline{\smash{\beta}}_{1}D$ and $\underline{\alpha}_{1}\underline{\alpha}_{2} = \underline{\smash{\beta}}_{1}\underline{\smash{\beta}}_{2}$. It follows that 
$\underline{\alpha}A = (\underline{\alpha}_{1}\underline{\alpha}_{2})(\underline{\smash{\gamma}}C) = \underline{\smash{\beta}}B = (\underline{\smash{\beta}}_{1}\underline{\smash{\beta}}_{2})(\underline{\delta}D)$, and so $\underline{\smash{\gamma}}C = \underline{\delta}D$.

Finally, we prove that $\underline{\smash{\gamma}}C = \underline{\delta}D$ is primitive.
Assume that $\underline{\smash{\gamma}}C = \underline{\delta}D$ is not a primitive solution, i.e., there exist vectors $\underline{\alpha}^{'}$, $\underline{\alpha}^{''}$, $\underline{\beta}^{'}$, $\underline{\smash{\beta}}^{''}$, $\underline{\smash{\gamma}}^{'}$, $\underline{\delta}^{'} \in \mathbb{Z}^{+}_{m}$ and multisets $C,D\in \mathbb{N}_{m}$ such that
$C = \underline{\alpha}^{'}C^{'}$, $D = \underline{\smash{\beta}}^{'}D^{'}$, $\underline{\smash{\gamma}} = \underline{\alpha}^{''}\,\underline{\smash{\gamma}}^{''}$,
$\underline{\delta} = \underline{\smash{\beta}}^{''}\underline{\delta}^{''}$,
$\underline{\alpha}^{'}\underline{\alpha}^{''} = \underline{\smash{\beta}}^{'}\underline{\smash{\beta}}^{''} \ne (1)$. 

First, let us suppose that $\underline{\alpha}^{'}\underline{\alpha}^{''} = \underline{\smash{\beta}}^{'}\underline{\smash{\beta}}^{''} = (t_{1}, t_{2}, \dots{})$, $t_{1} > 1$ and
$\underline{\alpha}^{'}_{1} = \underline{\alpha}^{(i_{0}+1)}_{1}\underline{\alpha}^{(i_{0}+2)}_{1} \dots{},$
$\underline{\alpha}^{'}_{2} = \underline{\alpha}^{(i_{0}+1)}_{2}\underline{\alpha}^{(i_{0}+2)}_{2} \dots{},$
$\underline{\smash{\beta}}^{'}_{1} = \underline{\smash{\beta}}^{(i_{0}+1)}_{1}\underline{\smash{\beta}}^{(i_{0}+2)}_{1} \dots{}$,
$\underline{\smash{\beta}}^{'}_{2} = \underline{\smash{\beta}}^{(i_{0}+1)}_{2}\underline{\smash{\beta}}^{(i_{0}+2)}_{2} \dots{}$.
It follows that
\[
\underline{\alpha}A = (\underline{\alpha}_{1}\underline{\alpha}_{2})(\underline{\smash{\gamma}}C) = \underline{\alpha}^{(1)}_{1}\underline{\alpha}^{(2)}_{1}\cdots{}
\underline{\alpha}^{(i_{0})}_{1}\underline{\alpha}^{'}_{1}\underline{\alpha}^{(1)}_{2}
\underline{\alpha}^{(2)}_{2}\cdots{} \underline{\alpha}^{(i_{0})}_{2}\underline{\alpha}^{'}_{2}
\underline{\alpha}^{''}\underline{\smash{\gamma}}^{''}\cdot \underline{\alpha}^{'}C^{'}
\]
\[
= (\underline{\alpha}^{(1)}_{1}\underline{\alpha}^{(2)}_{1}\cdots{}
\underline{\alpha}^{(i_{0})}_{1}\underline{\alpha}^{'})(\underline{\alpha}^{(1)}_{2}
\underline{\alpha}^{(2)}_{2}\cdots{} \underline{\alpha}^{(i_{0})}_{2}\underline{\alpha}^{''})
(\underline{\alpha}^{'}_{1}\underline{\alpha}^{'}_{2}\underline{\smash{\gamma}}^{''}C^{'})
\]
and similarly $\underline{\smash{\beta}}B = (\underline{\smash{\beta}}^{(1)}_{1}\underline{\smash{\beta}}^{(2)}_{1}\cdots{}
\underline{\smash{\beta}}^{(i_{0})}_{1}\underline{\smash{\beta}}^{'})(\underline{\smash{\beta}}^{(1)}_{2}
\underline{\smash{\beta}}^{(2)}_{2}\cdots{} \underline{\smash{\beta}}^{(i_{0})}_{2}\underline{\smash{\beta}}^{''})
(\underline{\smash{\beta}}^{'}_{1}\underline{\smash{\beta}}^{'}_{2}\underline{\delta}^{''}D^{'})$,
which is a contradiction because $i_{0}$ is maximal. 

One can prove a similar contradiction if $\underline{\alpha}^{'}\underline{\alpha}^{''} = \underline{\smash{\beta}}^{'}\underline{\smash{\beta}}^{''} = (1,n\dots{})$.

\section{Proof of Theorem \ref{theorem_3}}

$(i) \rightarrow (ii)$: Assume that $\log \lambda$, $\log \lambda^{'}$ are linearly independent over $\mathbb{Q}$, but \eqref{eq:4} is not a primitive solution.
Then then there exist
vectors $\underline{\smash{\gamma}}, \underline{\smash{\gamma}}^{'}, \underline{\smash{\mu}}, \underline{\smash{\mu}}^{'}, \underline{\nu},  \underline{\nu}^{'}$ and multisets $C, C^{'}\in \mathbb{N}_{m}$ such that 
\[
\lambda_{\underline{u},\underline{v}}^{(\underline{u})}S^{(\lambda,N)} = \underline{\smash{\mu}}C, 
\]
\[
\lambda_{\underline{u},\underline{v}}^{(\underline{v})} = \underline{\smash{\gamma}}\,\underline{\nu}
\]
\[
(\lambda^{'})^{(\underline{u}^{'})}_{\underline{u}^{'},\underline{v}^{'}}S^{(\lambda^{'},N^{'})} = \underline{\smash{\mu}}^{'}C^{'},
\]
%\[
%\underline{u},\underline{v}^{(\underline{v})} = %\underline{\smash{\gamma}}\,\underline{\nu},
%\]
\[
(\lambda^{'})^{(\underline{v}^{'})}_{\underline{u}^{'},\underline{v}^{'}} = \underline{\smash{\gamma}}^{'}\,\underline{\nu}^{'},
\]
and  $\underline{\smash{\mu}}\,\underline{\nu} = \underline{\smash{\mu}}^{'}\,\underline{\nu}^{'} \neq (1)$.

Since $1 = R_{\underline{\smash{\mu}}C, \underline{\smash{\gamma}}\,\underline{\nu}} = R_{C, \underline{\smash{\gamma}}\,\underline{\smash{\mu}}\,\underline{\nu}}$ and so by Theorem 1.2 in [7], there exists $\lambda_{0} \ge 2$ such that $\underline{\smash{\gamma}}\,\underline{\smash{\mu}}\,\underline{\nu} = (1, \lambda_{0}, \dots{})$ where the coordinates are the powers of $\lambda_{0}$ and so the coordinates of $\underline{\smash{\gamma}}, \underline{\smash{\mu}}, \underline{\nu}$ are the powers of $\lambda_{0}$, too. Similarly, there exists $\lambda_{0}^{'} \ge 2$ such that $\underline{\smash{\gamma}}^{'}\,\underline{\smash{\mu}}^{'}\,\underline{\nu}^{'} = (1, \lambda_{0}^{'}, \dots{})$ where the coordinates are the powers of $\lambda_{0}^{'}$ and so the coordinates of $\underline{\smash{\gamma}}^{'}, \underline{\smash{\mu}}^{'}, \underline{\nu}^{'}$ are the powers of $\lambda_{0}^{'}$, too. Then with suitable $t, t^{'}\in \mathbb{Z}^{+}$, we have 
$\underline{\smash{\mu}}\,\underline{\nu} = (1, \lambda_{0}^{t}, \dots{})$ and $\underline{\smash{\mu}}^{'}\,\underline{\nu}^{'} = (1, (\lambda_{0}^{'})^{t^{'}}, \dots{})$
because $\underline{\smash{\gamma}} = (1, \dots{})$, $\underline{\smash{\gamma}}^{'} = (1, \dots{})$. Since $\underline{\smash{\mu}}\,\underline{\nu} =  \underline{\smash{\mu}}^{'}\,\underline{\nu}^{'}$, thus we have $\lambda_{0}^{'} = \lambda_{0}^{t/t^{'}}$. Since  $\lambda_{\underline{u},\underline{v}}^{(\underline{v})} = (1, \lambda^{u_{1}}, \dots{})$ and  $\underline{\smash{\gamma}}\,\underline{\smash{\nu}} = (1, \lambda_{0}^{s}, \dots{})$ for some $s \in \mathbb{Z}^{+}$, we have $\lambda = \lambda_{0}^{s/u_{1}}$. Similarly, since  $(\lambda^{'})^{(\underline{v}^{'})}_{\underline{u}^{'},\underline{v}^{'}} = (1, (\lambda^{'})^{u^{'}_{1}}, \dots{})$ and  $\underline{\smash{\gamma}}^{'}\,\underline{\smash{\mu}}^{'} = (1, (\lambda_{0}^{'})^{s^{'}}, \dots{})$ for some $s^{'} \in \mathbb{Z}^{+}$, we have $\lambda^{'} = (\lambda_{0}^{'})^{s^{'}/u_{1}^{'}} = \lambda_{0}^{ts^{'}/t^{'}u_{1}^{'}}$. Thus we have $\log \lambda = \frac{s}{u_{1}}\log\lambda_{0}$, $\log \lambda^{'} = \frac{s^{'}t}{t^{'}u_{1}^{'}}\log\lambda_{0}$ which is a contradiction because $\log \lambda$, $\log \lambda^{'}$ are linearly independent over $\mathbb{Q}$.

$(ii) \rightarrow (iii)$: Assume that \eqref{eq:4} is a primitive solution, but it follows from the first construction. Then there exist
vectors $\underline{\smash{\mu}}, \underline{\smash{\mu}}^{'}\in \mathbb{Z}_{m}^{+}$ and multiset $C\in \mathbb{N}_{m}$ such that $\lambda_{\underline{u},\underline{v}}^{(\underline{u})}S^{(\lambda,N)} = \underline{\smash{\mu}}C$,  $(\lambda^{'})^{(\underline{u}^{'})}_{\underline{u}^{'},
\underline{v}^{'}}S^{(\lambda^{'},N^{'})} = \underline{\smash{\mu}}^{'}C$ and $\underline{\smash{\mu}}{\lambda^{(\underline{v})}_{\underline{u},\underline{v}}} = \underline{\smash{\mu}}^{'}{(\lambda^{'})^{(\underline{v}^{'})}_{\underline{u}^{'},\underline{v}^{'}}} \neq (1)$. 
%where $\lambda_{\underline{u},\underline{v}}^{(\underline{u})}S^{(\lambda,N)} = %\underline{\smash{\mu}}^{'}%(\lambda^{'})^{(\underline{v}^{'})}_{\underline{u}^{'},\underline{v}^{'}}$. 

Since $\underline{\smash{\mu}}{\lambda^{(\underline{v})}_{\underline{u},\underline{v}}} = \underline{\smash{\mu}}^{'}{(\lambda^{'})^{(\underline{v}^{'})}_{\underline{u}^{'},\underline{v}^{'}}} \neq (1)$ and  $R_{\lambda^{(\underline{u})}_{\underline{u},\underline{v}}S^{(\lambda,N)},
\lambda^{(\underline{v})}_{\underline{u},\underline{v}}} = R_{(\lambda^{'})^{(\underline{u}^{'})}_{\underline{u}^{'},\underline{v}^{'}}S^{(\lambda^{'},N^{'})},\lambda^{(\underline{v}^{'})}_{\underline{u}^{'},\underline{v}^{'}}}$, it follows that \eqref{eq:4} is not a primitive solution which is a contradiction. 

$(iii) \rightarrow (i)$: Assume \eqref{eq:4} is not follow from the first type construction but $\log \lambda$, $\log \lambda^{'}$ are not linearly independent over $\mathbb{Q}$, i.e., there exist $s, s^{'}\in \mathbb{Z}^{+}$ with $\lambda^{s} = (\lambda^{'})^{s^{'}}$. 
 Since
\[
\lambda^{\frac{s}{\textnormal{ gcd}(s,s^{'})}} = (\lambda^{'})^{\frac{s^{'}}{\textnormal{ gcd}(s,s^{'})}},
\]
where $\textnormal{gcd}\left(\frac{s}{\textnormal{ gcd}(s,s^{'})}, \frac{s^{'}}{\textnormal{ gcd}(s,s^{'})}\right) = 1$, and so there exists $\lambda_{0} \ge 2$, $\lambda_{0}\in \mathbb{Z}^{+}$ such that $\lambda = \lambda_{0}^{\frac{s^{'}}{\textnormal{ gcd}(s,s^{'})}}$, $\lambda^{'} = \lambda_{0}^{\frac{s}{\textnormal{ gcd}(s,s^{'})}}$.    
%If $N\in \mathbb{Z}^{+}$, 
Then 
\[
S^{(\lambda,N)} = (1,\lambda^{N}, \lambda^{2N}, \dots{} ) \{0,1,\dots{} ,\lambda-1\} 
\]
\[
= (1,\lambda^{N}, \lambda^{2N}, \dots{} ) (1,\lambda_{0}, \lambda_{0}^{2}, \dots{} ,\lambda_{0}^{\frac{s^{'}}{\textnormal{ gcd}(s,s^{'})}-1}) \{0,1,\dots{} ,\lambda_{0}-1\}
\]
%If $N$ is infinite, then 
%\[
%S^{(\lambda,N)} = \{0,1,\dots{} ,\lambda-1\} 
%\]
%\[
%= (1,\lambda_{0}, \lambda_{0}^{2}, \dots{} ,\lambda_{0}^{\frac{s^{'}}{\textnormal{ gcd}(s,s^{'})}-1}) \{0,1,\dots{} ,\lambda_{0}-1\}. 
%\]
%Let $N^{'}\in \mathbb{Z}^{+}$, then 
and
\[
S^{(\lambda^{'},N^{'})} = (1,(\lambda^{'})^{N^{'}}, (\lambda^{'})^{2N^{'}}, \dots{} ) \{0,1,\dots{} ,\lambda^{'}-1\} 
\]
\[
= (1,(\lambda^{'})^{N^{'}}, (\lambda^{'})^{2N^{'}}, \dots{} ) (1,\lambda_{0}, \lambda_{0}^{2}, \dots{} ,\lambda_{0}^{\frac{s}{\textnormal{ gcd}(s,s^{'})}-1}) \{0,1,\dots{} ,\lambda_{0}-1\}.
\]
%If $N^{'}$ is infinite, then 
%%\[
%S^{(\lambda^{'},N^{'})} = \{0,1,\dots{} ,\lambda^{'}-1\} = (1,\lambda_{0}, \lambda_{0}^{2}, \dots{} \lambda_{0}^{\frac{s}{\textnormal{ gcd}(s,s^{'})}-1}) \{0,1,\dots{} ,\lambda_{0}-1\}. 
%\]
Then for some $\underline{\smash{\mu}}, \underline{\smash{\mu}}^{'}\in \mathbb{Z}_{m}$ 
\[
\lambda^{(\underline{u})}_{\underline{u},\underline{v}}S^{(\lambda,N)} = \underline{\smash{\mu}}\{0,1,\dots{} ,\lambda_{0}-1\}, \hspace*{5mm} (\lambda^{'})^{(\underline{u}^{'})}_{\underline{u}^{'},\underline{v}^{'}}S^{(\lambda^{'},N^{'})} = \underline{\smash{\mu}}^{'}\{0,1,\dots{} ,\lambda_{0}-1\}
\]
which imply that
\[
 \mathbb{N} = \lambda^{(\underline{u})}_{\underline{u},\underline{v}}\lambda_{\underline{u},\underline{v}}^{(\underline{v})}S^{(\lambda,N)} = \underline{\smash{\mu}}\lambda_{\underline{u},\underline{v}}^{(\underline{v})}
 \{0,1,\dots{} ,\lambda_{0}-1\}
 \]
 \[
= (\lambda^{'})^{(\underline{u}^{'})}_{\underline{u}^{'},\underline{v}^{'}}\lambda^{(\underline{v}^{'})}_{\underline{u}^{'}, \underline{v}^{'}}S^{(\lambda^{'},N^{'})} = \underline{\smash{\mu}}^{'}\lambda^{(\underline{v}^{'})}_{\underline{u}^{'}, \underline{v}^{'}}\{0,1,\dots{} ,\lambda_{0}-1\}
\]
By Theorem \ref{theorem_1}, we have $\underline{\smash{\mu}}\lambda^{(\underline{v})}_{\underline{u},\underline{v}} = \underline{\smash{\mu}}^{'}\lambda^{(\underline{v}^{'})}_{\underline{u}^{'},
       \underline{v}^{'}}$, and so \eqref{eq:4} follows from the first type construction, which is impossible.

\end{document}